\documentclass[a4paper,11pt]{amsart}
\usepackage{amsfonts,amssymb,amsmath,amsthm,abstract,color}
\usepackage{mathrsfs}
\usepackage[ps,all,arc,rotate]{xy}
\usepackage[lmargin=1in,rmargin=1in,tmargin=1in,bmargin=1in]{geometry}
\usepackage{fancyhdr}
\usepackage{pb-diagram}
\usepackage{graphicx}
\usepackage{faktor}
\usepackage{enumitem}
\usepackage{xcolor}
\usepackage{mathtools}

\newtagform{brackets}{[}{]}
\usetagform{brackets}

\usepackage{hyperref}
\hypersetup{colorlinks=true,  colorlinks=false,       
	linkcolor=red,          
	citecolor=red,
	citebordercolor=blue,
	urlcolor=magenta} 

\usepackage{titlesec}
\usepackage{lmodern}
\usepackage{pb-diagram}
\usepackage{pgf}
\usepackage{epigraph}
\usepackage{graphicx}
\usepackage{tikz}
\usepackage{tikz-cd}
\usepackage{comment}
\usepackage{verbatim}
\usepackage{mathtools}
\usepackage{titlesec}
\usepackage{chngcntr}
\usepackage{faktor}

\newcommand{\ra}{\rightarrow}      

\renewcommand{\geq}{\geqslant}
\renewcommand{\leq}{\leqslant}

\def\cB{\mathcal B}\def\cC{\mathcal C}
\def\cF{\mathcal F}

\def\cO{\mathcal O}

\def\AA{\mathbb A}\def\CC{\mathbb C}

\def\NN{\mathbb N}\def\PP{\mathbb P}
\def\QQ{\mathbb Q}\def\RR{\mathbb R}

\def\ZZ{\mathbb Z}

 \def\U{\mathrm{U}}\def\GL{\mathrm{GL}} \def\SL{\mathrm{SL}}
\def\PGL{\mathrm{PGL}}

\newcommand{\ten}{\otimes}

\newcommand{\ngit}{\! \mathbin{\text{\rotatebox[origin=c]{70}{\scalebox{1.2}{$\approx$}}}} \!}

\DeclareMathOperator{\Bl}{Bl}
\DeclareMathOperator{\Proj}{Proj}
\DeclareMathOperator{\quot}{Quot}

\DeclareMathOperator{\Lie}{Lie}

\DeclareMathOperator{\diag}{diag}
\DeclareMathOperator{\proj}{Proj}
\DeclareMathOperator{\Hom}{Hom}

\DeclareMathOperator{\sym}{Sym}

\DeclareMathOperator{\stab}{Stab}

\DeclareMathOperator{\wt}{wt}

\DeclareMathOperator{\id}{id}

\DeclareMathOperator{\Aut}{Aut}
\DeclareMathOperator{\Nil}{Nil}

\newcommand{\gr}{\mathrm{gr}}

\usepackage{endnotes}
\usepackage{textcomp}

\newcounter{thm}[section]

\titleformat{\section}[block]
{\centering \normalfont\scshape}{\thesection .}{1em}{}

\newtheorem{thm}[subsubsection]{Theorem}
\counterwithin*{subsubsection}{subsection}
\counterwithin*{subsubsection}{section}

\newtheorem{defn}[subsubsection]{Definition}
\newtheorem{prop}[subsubsection]{Proposition}
\newtheorem{cor}[subsubsection]{Corollary}
\newtheorem{lem}[subsubsection]{Lemma}
\newtheorem{conj}[subsubsection]{Conjecture}
\newtheorem{example}[subsubsection]{Example}
\newtheorem{rmk}[subsubsection]{Remark}
\newtheorem{assertion}{ \color{red} Assertion}
\newtheorem{idea}[subsubsection]{Idea}
\newtheorem{aim}[subsubsection]{Aim}
\newtheorem*{qn}{ \color{blue} Question}
\newtheorem{cl}[subsubsection]{Claim}
\newtheorem{ntt}[subsubsection]{Notation}
\newtheorem{nt}[subsubsection]{Note}

\newtheorem{slg}{Slogan}
\newtheorem{exe}[subsubsection]{Exercise}
\newtheorem{ass}[subsubsection]{Condition}

\newcommand{\bass}{\begin{ass}}
	\newcommand{\eass}{\end{ass}}

\newenvironment{reference}{\paragraph{Reference:}}{\hfill$\square$}

\newcommand{\bexe}{\begin{exe}}
	\newcommand{\eexe}{\end{exe}}

\newcommand{\brf}{\begin{reference}}
	\newcommand{\erf}{\end{reference}}
\newcommand{\bnt}{\begin{nt} \normalfont}
	\newcommand{\ent}{\end{nt}}
\newcommand{\bntt}{\begin{ntt}}
	\newcommand{\entt}{\end{ntt}}
\newcommand{\bcl}{\begin{cl} \normalfont }
	\newcommand{\ecl}{\end{cl}}
\newcommand{\bqn}{\begin{qn} \normalfont  \begin{bf} }
		\newcommand{\eqn}{\end{bf} \end{qn}}
\newcommand{\bid}{\begin{idea}}
	\newcommand{\eid}{\end{idea}}

\newcommand{\bas}{\begin{assertion} \normalfont \begin{bf} }
		\newcommand{\eas}{\end{bf} \end{assertion}}
\newcommand{\bcr}{\begin{cor}}
	\newcommand{\ecr}{\end{cor}}
\newcommand{\bex}{\begin{example} \normalfont}
	\newcommand{\eex}{\end{example}}
\newcommand{\blm}{\begin{lem}}
	\newcommand{\elm}{\end{lem}}
\newcommand{\bthm}{\begin{thm}}
	\newcommand{\ethm}{\end{thm}}
\newcommand{\bcd}{\begin{tikzcd}}
\newcommand{\ecd}{\end{tikzcd}}
\newcommand{\bdf}{\begin{defn} \normalfont}
	\newcommand{\edf}{\end{defn}}
\newcommand{\bpp}{\begin{prop}}
	\newcommand{\bam}{\begin{aim}}
		\newcommand{\eam}{\end{aim}}
	
	\newcommand{\rar}{\rightarrow}
	\newcommand{\brm}{\begin{rmk} \normalfont }
		\newcommand{\erm}{\end{rmk}}
	\newcommand{\epp}{\end{prop}}
\newcommand{\bpf}{\begin{proof}}
	
	\newcommand{\hrar}{\hookrightarrow}
	
	\newcommand{\epf}{\end{proof}}

\newcommand{\kom}{$\mathcal{\kom}$}

\newcommand{\pno}{\mathbb{P}^{1}}

\newcommand{\overbar}[1]{\mkern 1.5mu\overline{\mkern-1.5mu#1\mkern-1.5mu}\mkern 1.5mu}

\newcommand{\bnu}{\begin{enumerate}}
	\newcommand{\enu}{\end{enumerate}}

\newcommand{\Gm}{\mathbb{G}_{m}}
\newcommand{\Ga}{\mathbb{G}_{a}}

\newcommand{\bpm}{\begin{pmatrix}}
	\newcommand{\epm}{\end{pmatrix}}
\newcommand{\eps}{\varepsilon}
\newcommand{\bfct}{\begin{fact}}
	\newcommand{\efct}{\end{fact}}
\newcommand{\bslg}{\begin{slg}}
	\newcommand{\eslg}{\end{slg}}
\newcommand{\Uh}{\hat{U}}
\newcommand{\git}{\mathbin{
		\mathchoice{\mkern-5mu/\mkern-6mu/\mkern-4mu}
		{\mkern-5mu/\mkern-6mu/\mkern-4mu}
		{\mkern-5mu/\mkern-6mu/\mkern-4mu}
		{\mkern-5mu/\mkern-6mu/\mkern-4mu}}}
\newcommand{\mc}{\mathscr}

\newcommand{\xddots}{%
	\raise 4pt \hbox {.}
	\mkern 6mu
	\raise 1pt \hbox {.}
	\mkern 6mu
	\raise -2pt \hbox {.}
}
\newcommand{\bgs}{\begin{gss}}
	\newcommand{\egs}{\end{gss}}
\def\acts{\curvearrowright}

\newcommand{\bit}{\begin{itemize}}
	\newcommand{\eit}{\end{itemize}}

\DeclareMathOperator{\Quot}{Quot}
\DeclareMathOperator{\End}{End}
\DeclareMathOperator{\Ext}{Ext}

\def\b{\beta}

\def\l{\lambda}
\def\m{\mu}

\newcommand{\Zb}{Z_\beta}
\newcommand{\Yb}{Y_\beta}
\newcommand{\Sb}{S_\beta}

\newcommand{\Pb}{P_\beta}
\newcommand{\Ybss}{Y_\b^{ss}}

\renewcommand{\Gamma}{\varGamma}
\renewcommand{\phi}{\varphi}

\mathchardef\ordinarycolon\mathcode`\:
\mathcode`\:=\string"8000
\begingroup \catcode`\:=\active
\gdef:{\mathrel{\mathop\ordinarycolon}}
\endgroup

\newcommand{\bcon}{\begin{conj}}
	\newcommand{\econ}{\end{conj}}

\newcommand{\HN}{Harder-Narasimhan }

\titleformat{\subsection}[runin]
{\normalfont\bfseries}{\thesubsection .}{1em}{}

\titleformat{\subsubsection}[runin]
{\normalfont\bfseries}{\thesubsubsection .}{1em}{}

\swapnumbers

\swapnumbers

\begin{document}
	
	\title[Moduli Spaces of Unstable Objects: Sheaves of Harder-Narasimhan length 2]{Moduli Spaces of Unstable Objects: \\ Sheaves of Harder-Narasimhan length 2}

\author{Joshua Jackson}
\address{J.\ Jackson  \\ Mathematics Department, Imperial College \\London \\ SW7 2AZ \\UK  \vspace{-5pt}}
\address{Heilbronn Institute for Mathematical Research\\ Bristol \\ UK}
\email{j.jackson@imperial.ac.uk}

	
	
	\maketitle

	\maketitle
	\begin{abstract}
		Given a moduli problem posed using Geometric Invariant Theory, one can use Non-Reductive Geometric Invariant Theory to quotient unstable HKKN strata and construct \lq moduli spaces of unstable objects', extending the usual moduli classifications. After giving a self-contained account of how to do this, we apply this method to construct moduli spaces for certain unstable coherent sheaves of HN length 2 on a projective scheme, which we call $\tau$-stable sheaves. This extends a previous result of Brambila-Paz and Mata-Guti\'{e}rrez for rank two vector bundles on a curve. 
	\end{abstract}

	\section{Introduction}

For more than fifty years, Geometric Invariant Theory (GIT) has been a central tool in moduli theory, and many important algebro-geometric moduli spaces can be constructed this way. Constructions generally all follow the same strategy. One starts by producing a scheme $X$ which parameterises the objects of interest, together with a redundancy given by the action of some \emph{reductive} linear algebraic group $G$. After describing the GIT-semistability condition $X^{ss} \subset X$ for this set-up in an object-intrinsic way, Mumford's GIT \cite{Mumford1994}, \cite{Newstead1978} then guarantees a projective quotient $X^{ss} \rar X\git G$ of the semistable locus, which is a good moduli space for these semistable objects. 
	
	One way of viewing such constructions, then, is as picking out a certain subclass of \lq semistable' objects whose members can be parameterised by a well-behaved space; the projectivity of the moduli space then corresponds to the assertion that this subclass is maximal, in the sense that no more members can be added to it. The remaining objects are called \emph{unstable},\footnote{The convention that \lq unstable' means \lq not semistable' rather than \lq not stable' is in some ways unfortunate, but by now far too well embedded to be avoided.} and, with few exceptions, moduli spaces parameterising them have yet to be constructed.     
	
	First, we provide a self-contained account of results of \cite{Hoskins} explaining how, subject to some fairly mild assumptions, moduli spaces for unstable objects can indeed be constructed. Secondly, to apply this approach to construct moduli spaces for unstable coherent sheaves of \HN length two, generalising a result of \cite{BrambilaPaz2013}. The main technical tools come from a new form of Geometric Invariant Theory for Non-Reductive Groups, developed in \cite{Doran2007} \cite{Berczi2020} \cite{Berczi2016} \cite{Berczi2017a} \cite{BDHK_handbook}\cite{BDK_grad_lin}\cite{Berczi2018a}

	\subsection{Moduli of unstable objects} \label{subsec intro moduli of unst obj}

	In the classical reductive GIT set-up described above, the moduli classification thus obtained says nothing at all about those objects whose corresponding points were in $X \setminus X^{ss}$, that is,  \emph{unstable}. While there are good geometric reasons not to try to include these objects in the same moduli space as the semistable ones, there is usually no good reason to ignore them altogether. We take the view that this is unsatisfactory, and that, rather than simply throw away the unstable points, the more natural thing to do is to regard the classical quotient $X\git G$ (the moduli space of semistable objects) as merely the first step in the solution of the moduli problem. 
	
	Just as classical GIT tells us which objects we can include in the moduli space of semistables, it also tells us how we should divide the problem further in order to make sense of the unstable part. This information is provided by the Hesselink-Kempf-Kirwan-Ness (HKKN) stratification, also called the instability stratification \cite{Hesselink1978}  \cite{Kempf1978} \cite{Kirwan1984} \cite{Ness1984}. As explained in \S \ref{subsec HKKN}, this stratification may be thought of as measuring how unstable points in $X\setminus X^{ss}$ are. Given a choice of invariant inner on product on $\Lie G$, one obtains a $G$-invariant locally closed stratification \[X = X^{ss}\sqcup \bigsqcup_{\beta \in \mc{B}} S_\b, \] which, roughly speaking, uses a normalised version of the Hilbert-Mumford function, building on work of Kempf \cite{Kempf1978}, to associate to each unstable point $x \in X\setminus X^{ss}$ a cocharacter of $G$ that is most responsible for the instability of $x$. In other words, one uses the fact that instability in GIT is structured phenomenon.  One can then try to find, for every instability index $\beta \in \mc{B}$, a quotient of the stratum $S_{\beta}$ by the $G$ action, and hence obtain a moduli space of objects with that instability index. A suitable slogan might be: we treat instability itself as a discrete invariant for the moduli problem. 
	
	The classical GIT quotients of $S_\b$ by $G$ are of course empty, at least with respect to the linearisation we started with, however due to work of Kirwan \cite{Kirwan1984} (see Theorem \ref{instab strat thm}), we may view the problem as equivalent to that of finding a quotient for the action of a certain parabolic group, $P_{\beta} \leq G$, on a subvariety $\Yb^{ss} \subset S_\beta$. Considering this latter question instead allows a greater range of linearisations to be brought into play, via twisting by characters. The trade-off for this is that \emph{$P_{\beta}$ is in general non-reductive}, meaning that Mumford's GIT cannot be used to construct a quotient. Hence, we will need to make use of the recent progress \cite{Berczi2020}  \cite{Berczi2017a} \cite{BDHK_handbook}\cite{BDK_grad_lin}\cite{Berczi2018a} \cite{Berczi2016},  \cite{Berczi2020} in producing a version of Mumford's GIT that works for non-reductive groups. 
	
	Non-reductive GIT encounters considerable \emph{prima facie} obstacles, the most famous of which is the non-finite generation of invariants seen in the Nagata counterexample \cite{Nagata1960} to Hilbert's 14th Problem. Nevertheless, the results of \cite{Berczi2020} provide a powerful machine that allows many of the good properties of the reductive case to be recovered, under mild technical assumptions. We will refer to these results collectively as the different versions of the $\Uh$ Theorem, summarised in \S \ref{sec Uhat}. Taking this as a starting point, in \cite{Hoskins} we prove the following. (See Theorem \ref{thm can qnt unstab strata} for a precise statement.)

	\bthm\cite{Hoskins} \label{intro thm can quotient unstable strata}   
	Let $X$ be a projective scheme, acted on by $G$ a reductive linear algebraic group. Suppose that the $G$ action is linearised with respect to an ample line bundle. Choose an invariant inner product $\|\cdot\|$ on $\Lie G$, and consider the associated HKKN stratification. 
	
	Then:
	\bnu
	\item  For each $\b \in \mc{B}$ such that a certain \lq semistability coincides with stability\rq\ condition holds for the unipotent radical of $P_\b$, \ref{*}, there is a projective good quotient by the action of the non-reductive group $P_\b$,
	
	\[\overbar{Y}_\b^{P_\beta-ss} \rar \overbar{Y}_\b \git_{\mc{L}_{\beta, \eps}} \Pb \] of a certain open locus in $\Ybss$, defined in the same way as in \cite{Berczi2020}.
	
	\item Further, this quotient map restricts to give a quasi-projective geometric quotient of the \emph{totally stable} locus $\Yb^{ts} \subseteq \overbar{Y}_\b^{P_\beta-ss}$, which we define at Definition \ref{defn Ys} and is given by a Hilbert-Mumford type criterion (\S \ref{HM subsec}). 
		
	\item Moreoever, if a certain \lq semistability coincides with stability\rq\ condition  \ref{ss=s R} holds for the reductive part of $P_\beta$, then $\overbar{Y}_\b^{P_\beta-ss}=\Yb^{ts}$, so that this geometric quotient of the totally stable locus is projective.

\enu 	
	\ethm 
	
Thus, one might summarise the strategy of moduli of unstable objects in the following way:
	
	\begin{enumerate}[label=\Roman*)]
		\item Pose the given moduli problem as a GIT problem.\footnote{Of course, the sets of objects in which we are interested will usually come in unbounded families, and in any given GIT-setup we will only see finitely many strata. So strictly speaking what we really need is an infinite family of GIT problems that somehow approximate our moduli problems asymptotically, so that we can work with an asymptotic instability stratification, in the manner of \cite{Hoskins2018}.}
		\item Interpret the instability stratification associated to this GIT problem in an intrinsic way. That is: convert the purely geometric-invariant-theoretic notion of HKKN type into a moduli-theoretic instability type. This gives a discrete classification, analogous to fixing the Hilbert polynomial or similar invariant.
		\item Use non-reductive GIT as set out in \cite{Berczi2020} to quotient the unstable strata, obtaining moduli spaces of objects of fixed instability type.

	\end{enumerate}
	
	The output of this is a moduli problem that has been split into (usually an infinite number of) manageable pieces, each with its own moduli space. Each isomorphism (or equivalence) class of object should appear in one, and only one, of these moduli spaces.

	\subsection{Unstable coherent sheaves}
	
	The remainder of this paper deals with our main application. Let $(X,\mc{O}_X(1))$ be a projective scheme with an ample line bundle, and suppose we wish to study the moduli of coherent sheaves on $X$. For sheaves that are \emph{semistable} (Definition \ref{ss sheaf defn}), this is a very well studied question: \cite{Huybrechts2010} is an excellent introduction to the vast literature on this topic. Unstable sheaves, by contrast, have received almost no attention: to our knowledge the only result is \cite{BrambilaPaz2013}. The latter result constructs a moduli space for rank two unstable vector bundles on a smooth projective curve, and proves a Torelli-type theorem in this context. The existence of this moduli space may be recovered from Theorem \ref{intro main thm} below. Our aim here is to try to understand unstable sheaves in a systematic way, using the strategy given above. 
	
	The now-standard GIT construction of the moduli space of semistable sheaves is due to Simpson \cite{Simpson1994}, and it is this construction that we take as our starting point, giving a brief recap in \S \ref{moduli ss sheaves}. The HKKN stratification coming from Simpson's construction is also well understood. As outlined in \S \ref{sheaves strats}, Nitsure and Shatz have shown \cite{Nitsure2011} that there is a second stratification on the parameter space, by Harder-Narasimhan (HN) type (Definition \ref{coprime HN defn}). As both the HN and HKKN stratifications both in some sense measure the failure of a sheaf to be semistable, it is reasonable to suppose that they might coincide, and due to work of Hoskins and Kirwan \cite{Hoskins2012} \cite{Hoskins2018}, we know that this is essentially the case  (Theorem \ref{thm comp strat}). Thus, in the case of unstable sheaves, parts (I) and (II) of the above strategy have already been completed, and we give a partial solution to (III), using Theorem \ref{intro thm can quotient unstable strata}.
	
	However, the quotient problem (III) for sheaves cannot be solved by a straightforward application of Thereoem \ref{intro thm can quotient unstable strata}. Indeed, we show in \S \ref{subsec automorphism gps of unstable sheaves} that the condition \ref{*} essentially never holds for unstable sheaves, so that instead we must use a different version of the $\Uh$ Theorem: following \cite{Berczi2020} we employ a kind of non-reductive \lq partial desingularisation' process (the analogy is with \cite{Kirwan1985}) to achieve a state of affairs in which the quotient may be performed. The output of this blow-up process may be thought of as a refinement of the HKKN stratification in the manner of \cite{Berczi2018}, with moduli spaces constructed for each refined HKKN-type. Thus, in the case of sheaves, one obtains some refinement of the HN type, which must be interpreted sheaf-theoretically. For sheaves of HN length two, we show in \S \ref{main section} that this refined invariant corresponds to fixing the dimension of the endomorphism group of the sheaf. 
	
	Whilst one may in principle perform this non-reductive partial desingularisation process for sheaves of arbitrary HN length, the sheaf-theoretic interpretation of the resultant refined HKKN stratification as a refinement of HN type does not appear to be tractable for sheaves of \HN lengths greater than two, for technical reasons explained at \S \ref{blow up discussion} and in Remark \ref{rmk higher length stab interp}. In addition, it is only for sheaves of \HN length 2 that the totally stable locus $Y^{ts}_\beta$ is non-empty. Thus, we restrict our attention to these sheaves, and leave the general case to be dealt with in a later work \cite{BHJK}, via a different non-reductive GIT construction. 
	
After restricting to sheaves of HN length 2 with fixed HN type $\tau$, we introduce the condition of \emph{$\tau$-stability} (Definition \ref{defn tau-stable sheaf}), which captures the NRGIT property of being totally stable in a sheaf-theoretic way, just as (semi)stability for sheaves captures GIT (semi)stability. For HN types of length 2 that are coprime, in the sense of Definition \ref{coprime HN defn}, $\tau$-stability is equivalent to indecomposability. Thus, our main result is:
	
	\bthm \label{intro main thm} Let $(X,\mc{O}_X(1))$ be a projective scheme with an ample line bundle, and let $\tau =(P_1,P_2)$ be a Harder-Narasimhan type of length 2. Then 
	\bnu 
	\item There is a quasi-projective moduli space $M_{\tau,d}$ for those $\tau$-stable (Definition \ref{defn tau-stable sheaf}) coherent sheaves such that $\dim\End(E) = d$. 
\item This moduli space has a canonical projective completion $M_{\tau,d} \subset \overbar{M}_{\tau,d}$, constructed via the nonreductive partial desingularisation procedure of \S \ref{blow up discussion}, such that the map \[M_{\tau,d} \rar M^{ss}_{P_1}\times M^{ss}_{P_2},\] to the product of two moduli spaces of semistable sheaves, which takes a sheaf to the associated graded of its HN filtration, extends to $\overbar{M}_{\tau,d}$.
\item Furthermore, if $\tau$ is coprime in the sense of Definition \ref{defn coprime}, then $M_{\tau,d}$ is the moduli space of indecomposable sheaves of type $\tau$ such that $\dim\End(E) = d$. 
\enu 
\ethm

	\subsection{Relation to previous work}
The basic idea of using non-reductive GIT to construct moduli spaces of unstable objects goes back at least as far as \cite{Hoskins2011} \cite{Hoskins2012}. However, the necessary results in Non-reductive GIT have only recently established. A thorough, if technical, account is given in \cite{Berczi2018}, where an inductive procedure is used to show that there exists a refinement of the HKKN stratification such that each refined stratum has a quasi-projective geometric quotient (that is, the quotient stack has a coarse moduli space). Our contribution is to give a closed-form expression for the open stratum of this construction, in the case covered by Theorem \ref{intro thm can quotient unstable strata}, and to interpret the GIT condition in a moduli-theoretic way, in the case of sheaves of HN length 2.

To our knowledge, the only previous construction of a moduli space for sheaves of fixed HN type is \cite{BrambilaPaz2013}, where the base scheme is a smooth projective curve, and the length two HN type is the type of a rank two vector bundle. The argument of \cite{BrambilaPaz2013} does not use GIT at all, either reductive or non-reductive, but rather constructs the vector bundle directly using a suitable sheaf of extensions. It is therefore of interest that the refinement of HN type thus obtained should coincide with ours.

As explained elsewhere in this paper, the case of HN length 2 enjoys certain special properties making it amenable to the methods we develop here, which do not hold in higher HN length. The construction of moduli spaces for unstable sheaves of HN length greater than two will be given in \cite{BHJK}, via a different NRGIT set-up.

		\subsection{Plan of the paper} 
		In \S \ref{git section} we recap the basics of reductive GIT and the HKKN stratification, concluding in \S\ref{hqotus} with a review of previous attempts to quotient unstable strata using reductive GIT . This motivates the use of Non-reductive GIT. In \S \ref{sec Uhat} we give an introduction to the results and terms from Non-reductive GIT that we will use later, including an extend discussion of the non-reductive partial desingularisation process, in \S \ref{blow up discussion}. In \S \ref{subsec nr qnts of unstable strata} we give a precise statement and proof of Theorem \ref{intro thm can quotient unstable strata}, and illustrate this with a simple example, in \S \ref{subsec pnts on pno}.
		
		We then turn to moduli of sheaves. We give the outline of Simpson's GIT construction for semistable sheaves in \S \ref{ss sh section}, and in \S \ref{sheaves strats} we review the work of Hoskins and Kirwan on comparing the HKKN and HN stratifications. This done, we conclude in \S \ref{main section} by giving an interpretation of the refined HKKN/HN stratification for length two sheaves, and proving Theorem \ref{intro main thm}.       
	~\\
	\begin{center} 
		\textbf{Acknowledgements} \end{center}  The author is grateful to Gergely B\'{e}rczi and Frances Kirwan for their very valuable help. Special thanks are due to Victoria Hoskins, whose thesis inspired this work, and for many helpful conversations.

	\subsection*{Notation} Throughout we work over an algebraically closed field $k$ of characteristic $0$. Unless otherwise stated, we will always use \lq sheaf' to mean coherent sheaf and \lq point' to mean closed point. Unless otherwise specified, following a helpful convention introduced by \cite{Doran2007}, $G$ will always be a reductive linear algebraic group, and $H$ an arbitrary (i.e.\ possibly non-reductive) linear algebraic group.

	\section{Reductive GIT and Instability Stratifications}\label{git section}

	\subsection{Recap of Reductive GIT} \label{recap of reductive GIT} Given a projective scheme $X$ equipped with the action of a linear algebraic group, the question of whether, or in what sense, a quotient variety exists can be very subtle. If the group $G$ is reductive, the answer to this question was provided by Mumford in the 1960's \cite{Mumford1994}, and may be summarised as follows. More details may be found in the approachable introduction, \cite{Newstead1978}.
	
	Embedding $X \subseteq \PP^n$ using a very ample line bundle $L$, we choose a linearisation of the $G$ action, i.e.\ a lift of the action to the total space of $L$. This extends the action to $G \acts \PP^n$ via a representation $G\rar \GL_{n+1}$. We say the action is linear, and denote by  $\mc{L}$ the line bundle $L$ together with this choice of extra data.
	This induces an action on the sections of all non-negative tensor powers of $L$, and hence an action on the graded ring, \[G \acts \bigoplus_{i=0}^\infty H^0(X,\mc{L}^i),\] where $H^0(X,\mc{L}^i)$ is just $H^0(X,L^i)$ equipped with the structure of a $G$-representation. This in turn allows us to speak of the sub-algebra \[\bigoplus_{i=0}^\infty H^0(X,\mc{L}^i)^G \subseteq \bigoplus_{i=0}^\infty H^0(X,\mc{L}^i)\] of invariant sections of non-negative powers of $L$. Reductivity of $G$ ensures that this algebra is finitely generated, meaning that the projective spectrum of this algebra is a projective scheme. We then define the (reductive) GIT quotient,\[X\git_\mc{L} G \coloneqq \Proj \left(\bigoplus_{i=0}^\infty H^0(X,\mc{L}^i)^G\right).\] One obtains a rational map $X \dashrightarrow X\git_\mc{L}G$ induced via the $\Proj$ functor by the inclusion of the invariants into the whole algebra, whose domain of regularity is the open set $X^{ss} \subseteq X$, called the semistable locus\footnote{We will write e.g. $X^{G-ss}$ if we wish to make explicit the dependence on the group.}.  This rational map restricts to give a good, and in particular categorical, quotient $X^{ss} \rar X \git_\mc{L} G$, which, locally on affine opens, is given by taking Spec of invariants. It restricts further to give a geometric quotient of the stable locus, $X^s \subseteq X^{ss} \subset X$, which in good cases will coincide with the semistable locus.
	
	 As a topological space, we have \[X\git_{\mc{L}} G \cong X^{ss}/\sim_S\] where the relation, called $S$-equivalence, is given by \[x \sim_S y \iff \overbar{Gx} \cap \overbar{Gy} \cap X^{ss} \neq \emptyset.\] Hence to understand the quotient it suffices to understand the semistable locus. Calculating a generating set of invariants by hand, however, is seldom possible in practice. It is therefore crucial to the utility of this theory that one has the Hilbert-Mumford criterion, as in \cite{Mumford1994}, Chapter 2.1. This allows one to calculate $X^{ss}$, and hence understand the quotient, by understanding semistability with respect to cocharacters $\l :\Gm \rar G$, i.e.\ without having to calculate the subalgebra of invariants.
	
	Given a cocharacter $\l  : \Gm \rar G$ and $x\in X$, the valuative criterion of properness guarantees that the map $\l(t)\cdot x: \AA_t^1\setminus \{0\} \rar X$ can be extended uniquely over the origin. We denote the image of the origin by the limit\footnote{If $k = \CC$ this is the limit in the analytic sense.} $ \overbar{x} \coloneqq \lim_{t \rar 0}\l(t)x$. This point is a fixed point for the $\l(\Gm)$ action on $X$, and hence the action of $\l(\Gm)$ restricts to an action on the fibre $L_{\overbar{x}}$. 
	
	\bdf \label{dfn HM function} We define the value of the \emph{Hilbert-Mumford function} $\m(\l,x) \in \ZZ$ to be the weight of the $\Gm$ action on $L_{\overbar{x}}$ If we wish to make the dependence on the linearisation explicit, we will write this as  $\m_{\mc{L}}(\l,x)$. \edf
	
	Then the \emph{Hilbert-Mumford criterion} tells us that \[	X^{ss} = \{x \in X \mid \m(\l,x) \geq 0 \text{ for all cocharacters } \l : \Gm \rar G\}\]
	\[	X^{s} = \{x \in X \mid \m(\l,x) > 0 \text{ for all cocharacters } \l : \Gm \rar G\}. \]

	Rather than checking cocharacters of $G$ one at a time, we can work with a whole maximal torus $T \leq G$ at once. Let $\mathfrak{t}^+$ be a positive Weyl chamber of $\mathfrak{t}=\Lie(T)$. Choosing coordinates $x_0,\dots x_n$ on $\mathbb{P}^n$ that diagonalise the $T$ action, the induced representation $T\rar \GL_{n+1}$ is determined by the weights $w_i \coloneqq \text{wt}_T(x_i) \in \mathfrak{t}^*$,  for $i=0,\dots, n$. 
	
	\bdf We associate to $x = [a_0:\dots a_n]\in X \subseteq \PP^n$ the convex polytope\[ \text{conv}_T(x) \coloneqq \text{conv}(w_i \mid a_i \neq 0) \subset \mathfrak{t}^* \] i.e.\ the convex hull of the set of weights $w_i$ whose corresponding coordinates $a_i$ are non-zero. \edf 
	
	It is not hard to show (\cite{Dolgachev2003}) that for tori the Hilbert-Mumford criterion of semistability can be reformulated to say: \[x \in X^{T-ss} \iff 0 \in \text{conv}_T(x),\]
	\[x \in X^{T-s} \iff 0 \in \text{Interior}(\text{conv}_T(x)).\]
	and that a point $x$ is (semi)stable for $G$ if and only if all points of its $G$-orbit are (semi)stable for $T$. Taken together, all of this reduces the calculation of $G$-(semi)stability to a finite problem, and gives a hint as to how we might measure instability.
	
	\brm A point that is not usually emphasised, but will be important later, is that the sense of \lq interior' meant above is \emph{with respect to the ambient space, $\mathfrak{t}^*$}. That is to say, if the polytope is not full-dimensional, then it is considered to have empty interior.  
	
	\erm

	\subsection{HKKN stratification} 
	\label{subsec HKKN}
	In fact, we can obtain more from reductive GIT than just a moduli space of semistables, via the Hesselink-Kempf-Kirwan-Ness (HKKN) stratification, or instability stratification, which we now recall. This stratification was used in \cite{Kirwan1984}, where it was used to calculate the cohomology of reductive GIT quotients (see e.g.\ \cite{Kirwan1989}) Our use for the stratification, however, is different: we regard it as a way to break the \lq full moduli problem' into manageable pieces. 
	
	Take $G \acts X \subseteq \PP^n$, a reductive group acting linearly on a projective variety, and notation as in the previous section. The stratification also depends on piece of additional data: we fix a norm $\|-\|$ on the set of cocharacters of $G$. That is, we take a norm on $\Lie G$ that is invariant under the action of the Weyl group.\footnote{In our main application, as well as many others, $G=\SL_r$, so the Killing form is the only such choice up to scaling.}

	\brm If $X$ is smooth, the HKKN stratification also admits a symplectic description in terms of the Morse stratification associated to the norm square of a suitable moment map. We will not use this description here. Details may be found in \cite{Kirwan1984}. \erm
	
	\bdf The index set $\mathscr{B}$ for the HKKN stratification is the set of points $\b \in \mathfrak{t}^*$ that are the closest points to the origin of the convex hull of some subset of the $T$-weights $w_i$. \edf 
	
	The rough idea of the instability stratification is that it indexes points $x \in X \setminus X^{ss}$ according to how unstable they are, by associating a cocharacter that is in a sense most responsible for the instability of $x$. To make this idea precise: Kempf showed \cite{Kempf1978} that for each $x \in X$ there is a cocharacter that achieves the minimum of the normalised Hilbert-Mumford function \[M(\l,x) \coloneqq \frac{\m(\l,x)}{\|\l\|}\]
	
	\bdf We say that a cocharacter $\l :\Gm\rar G$ is \emph{optimal} for a point $x \in X$ if it achieves the minimum $M(-,x)$ over all cocharacters of $G$. \edf

	Kempf shows that optimal cocharacters are unique up to conjugacy by a certain parabolic group $P_x$ associated to $x$, and each has a unique representative in $\mathfrak{t}^+$, a fixed positive Weyl chamber. This parabolic is given by \[P_x = P_\l \coloneqq \{g \in G \mid \lim_{t\rar 0}\l^{-1}(t)\cdot g\cdot\l(t)\text{ exists in $G$}\}\] where $\l$ is any optimal cocharacter.

	\brm 
	\label{rmk filtration parabolics} In the case we consider later, where $G = \SL(V)$, a cocharacter $\l: \Gm\rar G$ gives a splitting $V=\oplus_{r\in\ZZ}V_r$ of $V$ into weight spaces for the $\Gm$ action. Choosing a basis such that $\l(t) = \diag(t^{a_0},\dots,t^{a_n})$ with $a_i>a_{i+1}$, the parabolic $P_\l$ is the group of block upper triangular matrices preserving the filtration of $V$ whose pieces are $V^{i
	}\coloneqq \oplus_{r \geq i}V_r$. More generally, if $G \leq \GL_n$ for some $n$, then $P_\l$ is the intersection with $G$ of the standard parabolic of $\GL_n$ consisting of elements preserving the weight filtration of $\l$. 
	\label{parabolic remark} 
	
	\erm
	
	Next we need some notation that will be important later.

	\bdf \cite{Kirwan1984}
	
	\label{defns Zb Ybss etc}
		\bnu [label=\roman*)] 
	\item Using the inner product $\|.\|$, identify $\mathfrak{t}\cong \mathfrak{t}^*$. Then given $\b \in \mathfrak{t}^*$, we define \[\l_\b: \Gm \rar G\] to be the cocharacter corresponding $\beta$ under this identification. 
 \item We define $\Zb$ to be the closed subscheme of the fixed point locus $X^{\b(\Gm)}$ consisting of the points for which $\m(x,\l_\b) = -\|\b\|^2$. Considering $\b \in \mathfrak{t}$, we see that $\Zb$ is those points of $X$ whose $T$-weights lie on the hyperplane $H_\b$ perpendicular to $\b$ with respect to our chosen inner product.
	\item We define $p_\b :X\rar X^{\b(\Gm)}$ to be the retraction  \[p_\b(x) \coloneqq \lim_{t\rar 0} \l_\b(t)\cdot x.\] 
	
	\item We define the locally closed subscheme  $\Yb= p_\b^{-1}(\Zb)$. Equivalently, $\Yb$ is the set of points in $X$ whose $T$-weights all lie on the half-space \[H_\b^+ = \{w \in \mathfrak{t} \mid \langle v,\b \rangle \geq \|\b\|^2\}\] and at least one of whose weight lies on $H_\b$.
	
	\item We further define $$\Yb^{ss} = \{ x \in \Yb \mid \lambda_{\beta} \text{ is optimal for } x \}, $$  $$\Zb^{ss}  = \{x \in \Zb \mid \lambda_{\beta} \text{ is optimal for } x \},$$ $$S_{\beta} = G\cdot \Yb^{ss}.$$
	
	\item Finally we define $\overbar{\Yb}$ to be the closure of $\Yb$ (or equivalently $\Ybss$) in $X$.
	
	\enu 
	\edf 
	
	\bdf 
	\label{canonical lin defn}
	Denote by $\stab\beta$ the stabiliser of $\beta$ under the coadjoint action of $G$, and note that $\Zb$ is invariant under the action of $\stab\b$. We restrict the linearisation to the action of $\stab \beta$ on $\Zb$, and twist it by the character $-\beta$; this has the effect of shifting the weights by this vector. We denote this linearisation by $\mc{L}_\beta$, and we will call it the \emph{canonical linearisation}. 
	\edf
	
	\brm This character need not exist as a character of $G$, but the centrality of $\l_\b$ in $\stab\beta$ ensures that it exists as a character of $\stab\beta$, and, consequently, in $P_\beta$. This extra choice of character is one benefit of working with the non-reductive GIT problem we will consider later.\erm 
	
	\brm \label{Zbss remark}  The subscheme $\Zb^{ss}$ is exactly the semistable locus for the action of $\stab \beta$ on $\Zb$ with respect to the twisted linearisation, and $\Yb^{ss} = p_\beta^{-1}(\Zb^{ss})$
	(\cite{Kirwan1984}, 12.21). 
	\erm
	
	\bdf
	Given $\beta \in \mc{B}$, denote by $P_{\beta}$ the parabolic subgroup of $G$ associated to $\lambda_{\beta}$.
	\edf

	We can now summarise the relevant results of \cite{Kirwan1984} that give the stratification.
	
	\bthm \cite{Kirwan1984} \label{instab strat thm}
	Let $X$ be a projective scheme, acted on by a reductive linear algebraic group $G$, linearly with respect to an ample line bundle $\mc{L}$. Then, with notation as above, we have:
	\bnu
	\item There is a stratification of $X$ into disjoint locally closed $G$-invariant subschemes, $$X = \bigsqcup_{\beta \in \mc{B}} S_{\beta},$$  where the open stratum, $S_0$, is $X^{ss}$, and the closure of a stratum $S_{\beta}$ satisfies $$\overbar{S_{\beta}} \subset \bigcup_{\substack{\gamma \text{ s.t.} \\ \|\gamma\| \geq \|\beta\|}} S_{\gamma}.$$
	\item There is, for every $\beta \in \mc{B}$, a $G$-equivariant isomorphism \[\Yb^{ss} \times^{P_{\beta}}G \cong S_{\beta},\] where $\Yb^{ss}\times^{P_{\beta}}G$ is the quotient of $\Yb^{ss}\times G$ by the free action of $\Pb$ given by \[h\cdot(g,x) \coloneqq (gh^{-1},hx).\]
	\enu
	\ethm

	It is therefore natural to try to find, for every $\beta \in \mc{B}$, a quotient of the stratum $S_{\beta}$ by the $G$ action. By the second part of the theorem above, this is equivalent to the problem of finding a quotient for the action of $P_{\beta}$ on $\Yb^{ss}$. Thus, we are naturally led from a reductive GIT problem to a non-reductive one.

	\subsection{Reductive quotients of the unstable strata} \label{hqotus}

By way of motivation for the use of the $\Uh$ Theorem in this setting, we illustrate what happens if we try to quotient the unstable strata using reductive GIT. This was the approach taken in \cite{Hoskins2011}, before an adequate non-reductive GIT was available. Recall that in the best setting for GIT we have a reductive group $G$ acting on a projective scheme $X$ with ample linearisation $\mc{L}$, where the quotient morphism \[X^{ss} \rar X \git G\] is induced via the Proj functor by the inclusion of rings $$\bigoplus\limits_{r \geq 0}H^{0}(X,\mc{L}^{\otimes r})^G \hrar \bigoplus\limits_{r \geq 0}H^{0}(X,\mc{L}^{\otimes r}).$$ Outside of this setting, say when $L$ is not ample or $G$ is not reductive, we can still consider the morphism of schemes resulting from this inclusion, but this will not necessarily coincide with the usual GIT quotient in the sense of Mumford. Moroever, if the invariants are not finitely generated, the scheme obtained thereby will not be of finite type. A more tractable quotient may be obtained, following \cite{Hoskins2011}, by taking a finitely generated subalgebra \[R \subseteq \bigoplus\limits_{r \geq 0}H^{0}(X,\mc{L}^{\otimes r})^G\] and then defining the \lq generalised GIT quotient' to be the variety \[X\git_{\mc{L},R}G \coloneqq \proj R. \] This comes equipped with a natural rational map from the scheme $X$, defined by the corresponding inclusion of subrings.

Let $\overbar{Y}_\beta$ be the closure of $Y_\beta$ inside the ambient projective space, and define $$\hat{S}_\beta = G \times^{P_\beta} \overbar{Y}_\beta.$$ In \cite{Hoskins2011}, Hoskins constructs a linearisation on $\hat{S}_\beta$, which is also denoted by $\mc{L}_\beta$ and also called the canonical linearisation (cf.\ Definition \ref{canonical lin defn}). This has the property that there are isomorphisms of algebras

$$\bigoplus\limits_{n \geq 0}H^0(S_{\beta}, \mc{L}_{\beta}^{\otimes{n}})^G \cong \bigoplus\limits_{n \geq 0}H^0(\overbar{Y}_{\beta}, \mc{L}_{\beta}^{\otimes{n}})^{\Pb} \cong \bigoplus\limits_{n \geq 0}H^0(Z_{\beta}, \mc{L}_{\beta}^{\otimes{n}})^{\text{Stab}\beta}.$$

By the reductivity of Stab$\beta$, we know that the invariant sections on $\hat{S}_\beta$ are finitely generated. The canonical linearisation on $\hat{S}_{\beta}$ is therefore particularly amenable to the strategy of \lq generalised GIT\rq -- we can take $R$ to be, in each case, the whole subalgebra of invariants (Let us abuse notation and denote them all by $R$, since they are isomorphic.) On the other hand, this linearisation has the drawback that it is not usually ample.\footnote{See \cite{Hoskins2011}, Example 2.2.5} The result of this approach is summarised as follows.

\bthm (\cite{Hoskins2011}, $2.2$) 
\begin{enumerate}[label=\arabic*)]

	\item The semistable locus for the $\Pb$ action on $\overbar{Y}_\beta$ is $\Yb^{ss}$, and we have isomorphisms of generalised GIT quotients
	
	$$\hat{S}_\beta \git_{(\mc{L}_\beta, R)} G \cong \overbar{Y}_\beta \git_{(\mc{L}_\beta, R)} P_\beta \cong Z_\beta \git_{\mc{L}_\beta}\text{Stab}\beta.$$ 
	\item The quotient map $$\Yb^{ss} \rar \overbar{Y}_\beta \git_{(\mc{L}_\beta, R)} P_\beta $$ factors through the morphism $p_\beta: \Yb \rar \Zb$.
	\item The morphism $\Sb \rar \overbar{Y}_\beta \git_{(\mc{L}_\beta, R)} P_\beta $ is a categorical quotient for the action of $G$ on $S_\beta$.
	\enu 
	\ethm 
	
	The first point tell us the good news: in taking a quotient with respect to the canonical linearisation on $S_\beta$ we can essentially treat the problem as if it were classical GIT. That is, the quotient we are looking for is equivalent to a quotient of the action of the reductive group Stab$\beta$ on $\Zb$ with respect to the ample canonical linearisation. 
	
	However, the second point shows that each point of $\Yb^{ss}$ is identified in the quotient with its image under the retraction $p_\beta : \Yb \rar \Zb$. This is generally undesirable, and leads to GIT quotients being of much smaller dimension than they ought to be. For example, in the setting of moduli of unstable sheaves, this results in a sheaf being identified with the associated graded of its Harder-Narasimhan filtration; this means that the moduli spaces of unstable sheaves constructed in this way would be isomorphic to products of moduli spaces of semistable sheaves, with no information about the extensions that glue them together. This is plainly not satisfactory. A more straightforward illustration, where the moduli spaces resulting from this method of construction are just points, is given in \S \ref{subsec pnts on pno}, and we will see in \S \ref{sec moduli of unstable object} how we can use non-reductive GIT to improve on this.
	
	 In \S \ref{sec moduli of unstable object} we follow an idea due to Hoskins \cite{Hoskins2011}, first executed in \cite{Hoskins}. If we are to find a quotient morphism that does not factor through $p_{\beta}$, we must somehow make $\Zb$ unstable, otherwise we will have the problem that too many orbit closures meet there. Hence, a good approach would be to perturb the canonical linearisation, thus splitting each stratum into two pieces: a piece corresponding to $Z_{\beta}$, whose associated moduli space is the reductive quotient by $\stab \beta$ above, and another, open, piece\footnote{This will be $\Yb^{ts} = \Yb^s \setminus U_{\beta}Z_\beta^{s}$; compare to $X^0_{\min} \setminus U Z_{\min}$ in \S \ref{sec Uhat}.} which we will quotient using non-reductive GIT with respect to this perturbation of the canonical linearisation on $\overbar{Y}_{\beta}$. The problem arises, however, that perturbing the linearisation may potentially destroy the resemblance to a classical GIT problem noted in the theorem above, and creates a need for the methods of non-reductive GIT.

	\section{Non-reductive GIT: The $\Uh$ Theorem}\label{sec Uhat}
	
	Let us return to the situation above, in which a projective variety $X\subseteq \mathbb{P}^n$ is acted on linearly by a linear algebraic group $H$, now \emph{no longer assumed to be reductive}.

	\subsection{Non-Reductive GIT: The $\Uh$ Theorems}
	
 Suppose one attempts to generalise Mumford's theory in the most naive way, by proceeding exactly as in the reductive case: that is, declaring the quotient map to be the rational map from $X$ to the projective spectrum of the algebra of invariants, and taking $X^{ss}$ to be the locus of regularity of this map. Unfortunately, several things go wrong straight away. 
	
	\bit
	\item The subalgebra of invariants need not be finitely generated, so that this quotient may fail to be a scheme of finite type. This was first demonstrated by the Nagata counterexamples to Hilbert's 14th Problem, and simpler examples we later given by Roberts. \cite{Nagata1960} \cite{Roberts1990}
	\item Even when the subalgebra of invariants does happen to be finitely generated, there is no guarantee that the morphism as defined above have the topological properties we desire in a quotient map. In particular, the morphism can fail to be surjective. \cite{Doran2007}
	\item Unlike in the reductive case, taking invariants is not right exact. Consequently, non-reductive GIT quotients can in general behave badly with respect to $H$-equivariant closed embeddings.
	\item One way to interpret the Hilbert-Mumford criterion is as saying a reductive group $G$ has \lq enough' cocharacters to detect $G$-semistability. There is obviously no hope that an unmodified Hilbert-Mumford criterion can continue to hold for general non-reductive groups: for example, unipotent groups have no cocharacters at all.
	\eit 
	
	Nevertheless, all of these problems can be resolved under remarkably weak hypotheses. This generalisation of GIT relies on two crucial technical concepts, both of which we will show later are present in the cases of interest to us. We sketch it briefly here, and refer the reader in the first instance to \cite{Berczi2020}, as well as to \cite{Berczi2017a} \cite{Berczi2016}, for a more detailed account.
	
	\bdf \label{defn grading cocharacter} Let $H$ be an arbitrary linear algebraic group, and let $U$ denote its unipotent radical (i.e.\ its maximal connected normal unipotent subgroup). Since we are in characteristic zero, we can choose a reductive Levi factor $R \leq H$ and write $H=U\rtimes R$.  Let $\mathfrak{u} = \Lie(U)$. We say that a cocharacter $\l :\Gm \rar H$ \emph{grades $U$}, or that $U$ is \emph{positively graded} if it has trivial conjugation action on $H/U\cong R$ (i.e.\ its image lies in the centre of $R$) and the weights for its conjugation action on $\mathfrak{u}$ are all strictly positive.
	\edf
	
	\bdf  \label{defn adapted lin}
	We say that a linearisation is \emph{adapted} with respect to some choice of grading cocharacter $\l$, if the lowest weight of the $\l(\Gm)$ action on $X$ is strictly negative, and the rest are strictly positive.
	\edf
	
	\brm If the linearisation we are working with is not adapted, we can twist it by a suitable character to make it so. In the setting of moduli of unstable objects, we can always find a grading cocharacter, and adapted linearisation. Namely, if we wish to quotient $\Yb^{ss}$, we can use $\l_\b$ itself to grade the unipotent radical of $\Pb$, and twist the linearisation by $(1+\eps)(-\beta)$, where $0<\eps\ll 1$ and $\beta$ is the character dual to $\lambda_{\beta}$ under our choice of invariant inner product. \erm
	
	Before we can state the result from NRGIT that we will use, we need a few more pieces of notation, which we adopt in order to be consistent with \cite{Berczi2020}.
	
	\bdf
	Fix a grading $\l:\Gm\rar H$. Let $Z_{\min}$ be the closed subscheme of the fixed point locus $X^{\l(\Gm)}$ consisting of those points whose unique $\l$-weight is the minimal one. Then define $X^0_{\min} = p_\l^{-1}(Z_{\min})$, where $p_\l$ is as defined at Definition \ref{defns Zb Ybss etc}. 
	\edf
	
	We now quote the so-called $\Uh$ Theorem. The theorem get its name from a piece of notation: for a linear algebraic group $H$ with unipotent radical $U$ graded by a cocharacter $\l : \Gm \rar H$, we form the semidirect product $U\rtimes_\l\Gm$ and denote it by $\Uh$. There are various versions of the theorem that apply in different settings, but all of them require the following as a starting point. We will call these the $\Uh$ conditions.
	
	\bass \label{Uhat conds}
	Let $X$ be a projective scheme, with an action of a linear algebraic group $H$. Let $\l:\Gm\rar H$ be a grading cocharacter. Suppose that the action is linear with respect to an ample linearisation, and suppose that this linearisation is adapted in the sense of Definition \ref{defn adapted lin}. 
	\eass
	
	\bthm  \emph{($\Uh$ Theorem), \cite{Berczi2020} \label{Uhat thm}} Suppose the $\Uh$ conditions hold in the sense of \ref{Uhat conds}. Then with notation as in Definition \ref{defn grading cocharacter}, suppose further that \emph{semistability coincides with stability for the unipotent radical, $U$, of $H$}. That is, suppose we have\footnote{Requiring this condition for $X^0_{\min}$ is equivalent requiring it for  $Z_{\min}$: since $U$ is normalised by $\Gm$, the $U$-stabiliser dimensions cannot decrease along the retraction $p: X^0_{\min} \rar Z_{\min}$.}  
	\begin{equation} \label{*}  \tag{\lq\emph{ss=s for }$\Uh$\rq}
		\stab_U(x) = \{e\} \text{ for all } x \in Z_{\min}.
	\end{equation}
	
	Then
	
	\bnu 
	\item The $\Uh$-invariants are finitely generated, and the inclusion of the $\Uh$-invariant algebra induces a projective geometric quotient $$ X^{\Uh -s} \rar X \git_\mc{L} \Uh. $$ of an open subscheme of $X$ which we call the $\Uh$-stable locus: \[X^{\Uh-s} :=  X^{0}_{\min} \setminus UZ_{\min}.\]  
	\item Consequently the $H$-invariants are finitely generated, and the inclusion of the $H$-invariant subalgebra induces a good quotient of an open subscheme of $X$ (the $H$-semistable locus): $$X^{H-ss} \rar X \git_{\mc{L}} H,$$
	where $X \git_{\mc{L}} H$ is the GIT quotient of $X \git_{\mc{L}} \Uh$ by the induced action of the reductive group $H/\Uh \cong R/\Gm $ with respect to the induced linearisation. This good quotient of  $X^{H-ss}$ restricts to a geometric quotient \[X^{H-s} \rar X^{H-s}/H \subseteq X \git_{\mc{L}} H\] of the $H$-stable locus $X^{H-s}$.

	\enu 
	\ethm
	
	\brm The definitions of $X^{H-ss}$ and $X^{H-s}$ may be found in \cite{Berczi2020}; for our present purposes we need not determine these sets precisely, as it will be enough to work with a certain subset of $X^{H-s}$, which we explicitly describe in \S \ref{HM subsec}. \erm 
	
	For our main application, and others of interest, we will not have condition (\ref{*}). Fortunately, there are other versions of the $\Uh$ Theorem that can be used in such cases. These other versions follow a process that may be considered as analogous to the partial desingularisation procedure used in \cite{Kirwan1985} for reductive groups. An important caveat in this analogy is that, unlike in the reductive case, we must perform these blow-ups in order to obtain a quotient at all. 
	
	The version of the $\Uh$ Theorem that we will use in the construction of moduli spaces for sheaves of \HN length 2 is the following.
	
	\bthm \cite{Berczi2020} \label{Uhat for abelian U}
	Take the assumptions of \ref{Uhat conds}, without assuming (\ref{*}), and suppose that the unipotent radical $U$ is abelian. 
	
	Then there is an $H$-equivariant series of blow-ups of $X$, resulting in a variety $\pi:\widehat{X} \rar X$, with an $H$ action lifting the action on $X$, carrying an $H$-linearisation that is some small perturbation of $\pi^*\mc{L}$, such that \ref{Uhat conds} holds and in addition we have
	\begin{equation} \label{*_d}  \tag{$\star_d$}
	\exists d \in \NN \text{ such that } \dim\stab_U(x) = d \text{ for all } x \in \widehat{X}^0_{\min}.
	\end{equation}
	
	Furthermore, the conclusions (1) and (2) of \ref{Uhat thm} then hold for $\widehat{X}$.
	\ethm

	\brm \label{ rmk blow up non-ab U}
	As remarked in \cite{Berczi2020}, the above theorem can be combined with a certain quotienting-in-stages procedure, to yield the stronger version of the $\Uh$ Theorem that may be applied whenever we conditions of \ref{Uhat conds} are satisfied. For this we first must choose the a subnormal series \[0 \trianglelefteq U^1 \trianglelefteq \dots \trianglelefteq U^{\ell} = U, \] satisfying the condition that each successive quotient $U^{j+1}/U^j$ is abelian, and that the adjoint action of the distinguished cocharacter $\l(\Gm) \leq H$ on $\Lie U^{j+1}/U^j$ consists of a single weight space. Such a subnormal series can always be found: indeed, in the situation of moduli of unstable objects discussed in \S \ref{sec moduli of unstable object}, where we use the Killing form as our invariant inner product, we can use the derived series of $U$. One can then perform the quotient in stages, by first performing the blowups for the abelian group $U^1$ as in \ref{Uhat for abelian U}, and then following an iterative procedure for the remainder of $U$.
	\erm 
	
	The utility of \ref{Uhat for abelian U}, and its generalisation described above, lies in the fact that the map $\pi: \widehat{X} \rar X$ restricts away from the centres of the blow-ups to give an $H$-equivariant isomorphism of some nonempty open subset of $X$ with a nonempty open subset of $X$. Because being a good quotient is a local property, this set then has a good quotient by $H$. The difference in outcome between theorems \ref{Uhat thm} and \ref{Uhat for abelian U} is that in the former case we obtain an open subset of $X$ with a projective quotient, whereas in the latter case we obtain an open subset of $X$ whose quotient is only quasi-projective, but has a natural projective completion given by including orbits from $\widehat{X}$.

	\brm A final important property of this generalisation of reductive GIT, which we will not use here, is that Variation of GIT as in \cite{Dolgachev1998} \cite{Thaddeus1996}, can be shown to have an avatar in this setting \cite{Berczi2018a}. \erm

	\subsection{Hilbert-Mumford Criterion} \label{HM subsec} As we have observed, it is essential to have a tractable way of computing the stable and semistable loci if non-reductive GIT is to be useful in practice. Fortunately, subject to the conditions of the $\Uh$ Theorem, the Hilbert-Mumford criterion holds in exactly the sense we would like. A proof is given in full generality in \cite{Berczi2020}; here we will give a slightly simplified proof which is enough to establish the result in the case we will need it.

	Continuing to use the notation as in theorem \ref{Uhat thm}, write $\overbar{R} = R/\l(\Gm)$ for the quotient of the Levi by the grading cocharacter. Since this cocharacter acts trivially on $Z_{\min}$ by construction, we get an induced action $\overbar{R} \acts Z_{\min}$. Since the image of $\l : \Gm \rar H$ lies in the centre of $R$, the weight $w_{\min}$ with which $\l(\Gm)$ acts on $Z_{\min}$, \emph{a priori} only a character of $\l(\Gm)$, in fact corresponds to a character of $R$. Hence, twisting $\mc{L}$ by this character, we obtain a linearisation $\overbar{\mc{L}}$ of the action of $\overbar{R}$ on $Z_{\min}$. With respect to this linearisation, we consider the following condition: 
	
	\bdf \label{defn totally H-stable} We define the \emph{totally $H$-stable locus} to be \[X^{H-ts} := \{ x \in X^0_{\min} \setminus UZ_{\min} \mid \lim_{t \rar 0}\l(t)x  := p_{\lambda}(x) \in  Z_{\min}^{\overbar{R}-s}   \}\]
	\edf 
	
	Since $HZ_{\min} =UZ_{\min}$, this property and its negation are invariant under the action of $H$.
	
	\brm 
	We will see in \S \ref{main section} that this condition arises naturally when one considers unstable strata in the standard GIT construction of moduli of sheaves. Indeed, the HN-stable sheaves are exactly the sheaves satisfying this condition.
	\erm
	
	\bthm (Hilbert-Mumford Criterion) (cf. \cite{Berczi2020}, \cite{Hoskins}) \label{thm HM crit}
	
	Take the assumptions of Theorem \ref{Uhat thm}. Then
	
	\bnu \item All totally $H$-stable points are $H$-stable. That is, we have $X^{H-ts} \subseteq X^{H-s}$, and the totally $H$-stable locus therefore has a quasi-projective geometric quotient by $H$.

	\item If we have the following condition
	\begin{equation} \label{ss=s R} \tag{\lq \emph{ss=s for $\overbar{R}$}\rq} Z_{\min}^{\overbar{R}-s} = Z_{\min}^{\overbar{R}-ss} \end{equation} then there are equalities \[X^{H-ts} = X^{H-s} = X^{H-ss},\] so that the locus on the left has a projective geometric $H$-quotient.  
	
	\item Furthermore, (1) and (2) hold also in the situation of Theorem \ref{Uhat for abelian U}. That is to say, if we do not have condition (\ref{*}), but the unipotent radical is abelian, we obtain results as in (1) and (2) with $\widehat{X}$ replacing $X$. 
	
	\enu

	\ethm
	\bpf
	We show how the introduction of the assumption (\ref{ss=s R}) allows one to give an easier proof of this Hilbert-Mumford criterion, by simplifying the proof of Lemma 7.8 of \cite{Berczi2020}.  We give the proof for the situation of Theorem \ref{Uhat thm}; the proof for Theorem \ref{Uhat for abelian U} is the same,  following Remark 7.12 of \cite{Berczi2020}. 
	
	Adopting notation as in \cite{Berczi2020} \S7.4, where Proposition 7.4 gives a locally trivial $U$-quotient \[q_U:X^0_{\min} \rar X^0_{\min}/U.\]

	We choose  $s>0$ sufficiently large as determined by \cite{Berczi2020}, and define
	\[V = H^0(X,L^{\otimes s})^* \]
	\[V^U = (H^0(X,L^{\otimes s})^U)^*. \] We then obtain a locally closed immersion \[\iota : X^0_{\min}/U \hrar \PP(V^U).\] There is an induced $\l(\Gm)$ action on $\PP(V^U)$; we write $\PP(V^U)_{\min}$ for its minimal weight space and $\overbar{p}: \PP(V^U)_{\min}^0 \rar \PP(V^U)_{\min}$ for the induced retraction. 
	Now $\PP(V^U)$ has an action of the maximal torus $T\leq H$, which is reductive and hence enjoys a classical Hilbert-Mumford criterion. The idea of the proof in \cite{Berczi2020} is to bootstrap this Hilbert-Mumford criterion from $\PP(W^*)$ to $X^0_{\min}$ via $X^0_{\min}/U$. Recall from \S\ref{recap of reductive GIT} the formulation of the Hilbert-Mumford criterion which says that a point is stable for a torus action exactly when its weight polytope contains the origin in its interior.
	
	We must now compare, for a point $x \in X^0_{\min}$, the two weight polytopes $\text{conv}_T(x)$ and $\text{conv}_T(q_U(x))$. As is standard, we may replace $L$ with a positive tensor power, since the property that the weight polytope contain (or not contain) the origin is unaffected by this. Hence we may assume that the composition of the $U$-quotient map with $\iota$ corresponds to the inclusion of the $U$-invariant sections of $H^0(X,L^{\otimes s})$, and thus comes from a linear projection of the associated projective spaces $\PP(V) \dashrightarrow \PP(V^U)$. Thus, we have \[\text{conv}_T(q_U(x))  \subseteq \text{conv}_T(x),\] and \emph{a priori} this could be any subset.
	However, as observed in \cite{Berczi2020} the condition of $\l(\Gm)$ grading $U$ guarantees that $H^0(X,L^{\otimes s})_{\max} \subset H^0(X,L^{\otimes s})^U$ , so all weights of $\text{conv}_T(x)$ that are minimal for $\l$ are still present in $\text{conv}_T(q_U(x))$. Hence if $p(x) \in Z_{\min}^{\overbar{R}-s}$, then the same will be true after the $U$-quotient, i.e. $\overbar{p}\circ q_U(x) \in \PP(V^U)_{\min}^{\overbar{R-s}}$. Crucially, the weight polytope of the latter point will be of codimension one in $\mathfrak{t}^*$. Furthermore, we have $x \notin UZ_{\min}$, from which it follows that $q_U(x) \notin \PP(V^U)_{\min}$, since the $U$-quotient is geometric. Hence, $\text{conv}_T(x)$ contains at least one weight that is off the minimal weight space, and the same is true of $\text{conv}_T(q_U(x))$. The combination of these facts with well-adaptedness of the linearisation ensures that $q_U(x)$ satisfies the Hilbert-Mumford criterion for $T$-stability. Moreover, the same will be true for every point in the $H$-orbit of $x$, since having stable limit in $Z_{\min}$ is an $H$-invariant property. Hence, $x \in X^{H-s}$.

	\epf

	\subsection{The blow-up process in more detail} \label{blow up discussion}

	We have seen that, in quotienting by non-reductive groups, the dimensions of the stabilisers in the unipotent radical play an important role.
	For unstable sheaves of \HN length 2, we will interpret these stabilisers in \S \ref{subsec automorphism gps of unstable sheaves}, and see that we do not have condition \ref{ss=s R}. Consequently, we will need the full force of Theorem \ref{Uhat for abelian U}, and so we will need to understand the blow-up process used to prove that theorem. We adopt the assumptions of Theorem \ref{Uhat for abelian U}.
	
	In order for the blow-up process described in \cite{Berczi2020} to work, it is necessary that the centre of the blow-up is always of codimension at least two, or else the blow-up will be an equivariant isomorphism and the induction will stall. To prevent this problem from arising, here and in \cite{Hoskins}, we adopt a different approach, following an earlier arXiv version of \cite{Berczi2020}. In this approach, we blow up according to the dimension of $\Uh$-stabilisers rather than $U$-stabilisers - or equivalently, blowing up the locus in $UZ_{\min}$ of maximal $U$-stabiliser dimension. This implies that the centre of the blow-up is always of codimension at least two; that is unless we have $UZ_{\min} = X^0_{\min}$, in which case the reductive quotient $Z\git R$ is the appropriate $H$-quotient of $X$.
	
	Now let us examine this blow-up process. Proofs of the below statements, as far as they do not simply follow from the theorems quoted above, may be found in \cite{Berczi2020} and \cite{Hoskins}. 
	
\bdf  We will use the following notation for the blow-ups' centres. For $d\in \NN$, let \[C_{d}(X) = \{ x \in X^0_{\min} \mid \dim \stab_{\Uh}(x) = d\} \] and define \[d_{\max}(X)  = \max \{ \dim \stab_{\Uh}(x) \mid x \in X^0_{\min} \}.\] \edf  

With this notation, we perform a series of blow-ups to achieve condition (\ref{*_d}) for $U$- that is, to have a constant-dimensional $U$-stabiliser across all of $X^0_{\min}$. We begin by blowing up along the closure of the locus $C_{d_{\max}}(X)$. This gives us a new space $\pi^{(1)} : X^{(1)} \rar X$, with exceptional locus $E^{(1)}$. Because $\Uh$ is normal, the centre of the blow-up was $H$-invariant, so the new space has an $H$ action lifting the action on $X$, which we linearise by defining $$\mathscr{L}^{(1)} \coloneqq (\pi^{(1)})^* \mathscr{L}^{r}\ten \mathcal{O}_X(-E^{(1)}).$$ We choose $r \gg 0$, so that, since replacing a linearisation with a positive multiple makes no difference to (semi)stability, we may think of this linearisation as a perturbation of the original linearisation of the action on $X$ by a sufficiently small multiple of the exceptional divisor.

	We denote the new minimal weight space for the distinguished $\lambda :\Gm \rar H$ by $Z_{\min}(X^{(1)})$. Then let $$(X^{(1)})^0_{\min} =\{ x \in X^{(1)} \mid \lim_{t \rar 0} \lambda (t)x \in  Z_{\min}(X^{(1)})\}.$$ The arguments of \cite{Berczi2020} (see also \cite{Hoskins}\S 2.3.6-2.3.7) tell us that we now have $$ d^1_{\max}(X^{(1)}) \coloneqq \max \{ \dim \stab_{\Uh}(x)\mid x\in X^{(1)}\}  < d^1_{\max} .$$ So we repeat the process: we blow up $X^{(1)}$ along the closure of $C_{d_{\max}}(X^{(1)})$, obtaining $X^{(2)}$ with another action lifting that of $X$, which we again perturb by a small multiple of the exceptional divisor $E^{(2)}$. 
	
	Inductively, we see that this process will terminate with a space $ \widehat{X} :=X^{(l)}$ and a blow-down map $\pi: \widehat{X} \rar X$ such that the action of $U$ satisfies condition (\ref{*_d}). At this point we apply Theorem \ref{Uhat for abelian U}.

	\brm \label{blow up wts Gm}
	Let us pause at this point to illustrate the properties of the linearisation on the blow-up. For simplicity, assume that we have a $\Gm$ action on $\mathbb{P}^n$, with weights $w_i = \wt_{\Gm}(x_i)$, and that we are blowing up along a linear subspace $\mathbb{P}^{k} \subset \mathbb{P}^n$. Then the blow-up is \[\Bl_{\mathbb{P}^k} \mathbb{P}^n = \{([x_0:\dots:x_n],[\tilde{x}_{k+1}:\dots\tilde{x}_{n}] \in \mathbb{P}^n \times \mathbb{P}^{n-k-1} \mid x_i\tilde{x}_j = x_j\tilde{x}_i \}.\] With respect to the choice of linearisation given above, weights of the $\Gm$ action on the blow-up will be given by summing $r$ of the weights $w_0,\dots,w_n$ and then perturbing by one of the weights $w_{k+1},\dots w_{n}$: that is, we replace each weight $w_i$ by a small cluster of weights, and each weight in the cluster corresponds to a weight in the projectivised normal bundle of the centre of the blow-up. 
	
	In general, we will be blowing up a scheme $X$ that may be singular, non-reduced, etc, along a possibly singular locus. However, our choice of ample linearisation gives an embedding  $X \subset \mathbb{P}^n$, and we may perform the blow-ups along the relevant loci considered in this space, which, by the arguments of \cite{Berczi2020}, will be smooth because $\mathbb{P}^n$ is. This means that the identification of the exceptional divisor in $\Bl_{\mathbb{P}^k} \mathbb{P}^n$ with the projectivised normal bundle to the centre of the blow-up holds, and we may consider the exceptional divisor in the blow-up of $X$ to be the intersection of the proper transform of $X$ with this, which will therefore have weights as described above.
	
	\erm

	\brm For non-abelian $U$, as in the situation of Remark \ref{ rmk blow up non-ab U}, the idea is essentially the same. We first perform the blow-ups described above, with the abelian group $U^1$ in place of $U$. Having achieved condition (\ref{*_d}) for $U^1$, we then repeat the whole process for $U^{2}$, blowing up along closures of loci of the above form, and then inductively replacing $U^j$ with $U^{j-1}$ until we have a blow-up of $X$ satifying (\ref{*_d}) for the full unipotent radical $U$. \erm

	Of course, what one is really interested in is obtaining a quotient of an open subset of $X$. In the situation that will be of greatest concern to us, when we have posed a particular moduli problems using GIT, the $H$-orbits in $X$ will have some moduli-theoretic interpretation as equivalence classes of objects to be classified, whereas the points of $\widehat{X}$ may not readily admit such an interpretation. As a result, we will mainly use Theorem \ref{Uhat for abelian U} in the following form.
	
	\bcr \label{cr uh4} The blow down map $\pi : \widehat{X} \rar X$  restricts to the complement of the exceptional divisor $E$, to give an $H$-equivariant identification $X \setminus \pi(E) \cong \widehat{X} \setminus E$. Consequently, we may define in $X$ an open subset, $$X^{H-\widehat{ss}} \coloneqq \pi ( \widehat{X}^{H-ss} \setminus E ) \cong \widehat{X}^{H-ss} \setminus E. $$ This open set admits a quasi-projective good quotient, which we denote $$X^{H-\widehat{ss}} \rar X \ngit_{\mathscr{L}} H,$$ of which $\widehat{X}\git H $ is a projective completion. 
	
	\ecr
	
	Hence, the price that we must pay for the construction of moduli spaces when we do not have (\ref{*}) or (\ref{*_d}) for the unipotent radical is that the resultant spaces will not be compact; that is, they will be quasi-projective rather than projective, and the projective completions yielded by the above process will not necessarily readily admit modular interpretation.
	
	In order to determine $X^{H-\widehat{ss}}$, we must keep track of the minimal weight space and its basin of attraction at each stage of the process. To simplify notation slightly, let's suppose that we are at the very beginning of the process, i.e.\ that we are working with $X$ and have yet to do any blow-ups. We blow up along the closure of  $C_{d_{\max}}(X)$ to get $\pi^{(1)} : X^{(1)} \rar X$ with exceptional locus $E^{(1)}$. Now there are two cases.
	
	\begin{enumerate} [label=\arabic*)]
		\item If $Z_{\min} \not\subseteq C_{d_{\max}(X)}$, then the minimal weight space in the blow-up, $Z_{\min}(X^{(1)})$ is simply the proper transform of the old one, and $\pi^{(1)}$ gives an identification $$(X^{(1)})^0_{\min} \setminus E \cong \pi^{(1)}((X^{(1)})^0_{\min} \setminus E) = \{ x \in X^{0}_{\min} \mid \dim \stab_{U}^((p(x)) < d_{\max} \},$$ where we recall that $p : X^0_{\min} \rar Z_{\min}$ is the map given by $$x \longmapsto \lim_{t \rar 0}\lambda(t)x.$$
		\item If $Z_{\min} \subseteq C_{d_{\max}}(X)$, the situation is more complicated. Roughly speaking, since we are blowing up the whole minimal weight space, there will be nothing left of it, and consequently its place will be taken by that part of the exceptional divisor which lies over the old minimal weight space and corresponds to the lowest weight space that we haven't blown up entirely. More precisely, for $r \in \ZZ$, define $Z_{r}$ to be the $\lambda$ weight space of weight $r$, and then let $$W^0_r  = \{ x \in X^{0}_{\min} \mid \lim_{t \rar 0 }\lambda^{-1}(t)x \in Z_{r} \}. $$ That is to say, $W^0_r$ is those points of $X^0_{\min}$ whose highest weight is $r$. Now define  $$r_{\min} = \min \{ r \in \ZZ \mid W^0_r \not\subseteq C_{d_{\max}}(X) \}.$$ Then the new minimal weight space is that subset of $(\pi^{(1)})^{-1}(Z_{\min})$ given by $$p^{(1)} \circ (\pi^{(1)})^{-1} (W^0_r \setminus C_{d_{\max}}(X)), $$ where $p^{(1)}$ is the map induced by $p$ on the blow-up.
		\enu

		It is clear from the above that, although Corollary \ref{cr uh4} gives us an open subset of $X$ that has a categorical quotient, determining this locus in practice can be non-trivial - especially if $U$ is not abelian. Furthermore, if we are to find an interpretation of the categorical quotient $X^{H-\widehat{ss}} \rar X \ngit_{\mathscr{L}} H$ as a moduli space for certain objects in a moduli problem, we must understand the whole blow-up process in terms of the objects parameterised by $X$, and such conditions as being in a certain $C_{d}$ or $W^0_r$ should have some meaning intrinsic to the objects themselves.  
		
		Fortunately, as we will see in \S\ref{main section} in the blow-up process the sheaves of HN length $\ell =2$ one can show that the following additional condition holds. 
	\bdf \label{def p preserves U stabs} We say that\emph{$p$ preserves $U$-stabiliser dimension if}
		\begin{equation}
			\label{$p$ preserves U stab dim} \tag{$\dim\stab_U(x)=\dim\stab_U(p(x))$} \forall x\in X^0_{\min}\quad \dim\stab_U(x)  = \dim\stab_U(p(x)).
		\end{equation}
		
		\edf 
		
		The condition \ref{$p$ preserves U stab dim} simplifies the blow-up considerably, as one can easily show that when this condition holds the blow-up process never leaves case (1) above. This will be crucial to our application to sheaves.
		
		On the other hand, for $\ell > 2$ this condition does not hold, as observed in \cite{Hoskins}, and so we will sooner or later find ourselves in case (2). Interpreting the various loci of the form $C_{d}$ and $W^0_r$ that arise in the blow-up process for general length involves complicated nested Brill-Noether type conditions, and as such seems rather intractable. This, as well as the failure of condition \ref{ss=s R}, is the reason that we confine ourselves to the case of \HN length 2 in the present article. The general case will be dealt with in \cite{BHJK}.

		\brm \label{extrm} Although we will not consider the external case in detail here, one may use all of the $\Uh$ Theorems to perform the quotient of $X$ by $H$ using an external $\Gm$, in the following sense. We apply one of the  theorems above with $H$ replaced by $\hat{H} = H \rtimes \Gm$, and $X$ replaced by $X \times \pno$, linearising the action by tensoring $\mc{L}$ with $\mc{O}_{\mathbb{P}^1}(N)$ for $N \gg 0$. This yields a projective scheme, which if semistability equals stability for the unipotent radical is just $(X \times \pno) \git \hat{H}$, containing an open subscheme which is a geometric quotient of a certain $H$-invariant open subset $X^{H-\hat{s}} \subset X$ by $H$.   This approach is of course particularly useful if no admissible internal $\Gm$ exists, for example if $H$ is unipotent. \erm

		\section{Moduli Spaces of Unstable Objects} \label{sec moduli of unstable object}
		\label{sec quot HKKN strata}
	
		Having presented the essentials of Non-reductive GIT in the form of the $\Uh$ Theorem and its consequences, we turn to the question of constructing moduli spaces of unstable objects, outlined at \S\ref{subsec intro moduli of unst obj}. Recalling the set-up of \S\ref{subsec HKKN}, assume we have some moduli problem of interest posed using reductive GIT. That is to say, we have an action of a reductive group $G$ on a projective scheme $X$, such that orbits in $X$ correspond to isomorphism (or equivalence) classes of the objects we wish to classify. Assume also that the action is linearised with respect to some very ample line bundle, and that we have made a choice of invariant inner product $\|.\|$, so that we may consider the associated HKKN stratification of $X$. Then classical reductive GIT gives us a \lq moduli space of semistables' -a quotient of the open HKKN stratum- and our goal here is to find reasonable quotients for the unstable HKKN strata, and hence classify unstable objects too, by constructing for them their own moduli spaces. The $\Uh$ Theorem turns out to be ideally suited to this purpose.

			\subsection{Non-Reductive Quotients of unstable strata}
			
			\label{subsec nr qnts of unstable strata}

			Recall the notation of \S \ref{subsec HKKN}, and in particular of Definition \ref{defns Zb Ybss etc}. Since we are over an algebraically closed field of characteristic zero, we can write $P_\beta = U_\b \rtimes \stab\b$ where $U_\b$ is the unipotent radical of $P_\b$, and $\stab\b$ is as defined in \S \ref{subsec HKKN}. We will also need one further piece of notation.
			
			\bdf  \label{defn Ys}
			Recalling Remark \ref{Zbss remark}, we define $\Zb^{s}$ by analogy to be the stable locus for the action of $\stab \beta$ on $\Zb$, and $\Yb^{s} = p_\beta^{-1}(\Zb^{s})$. We then define by analogy with Definition \ref{defn totally H-stable} the \emph{totally stable locus}: \[Y_\b^{ts} := \Yb^s \setminus U_\b Z_\b^s.\]
			\edf 
			
			Let us now give a precise statement and proof of Theorem \ref{intro thm can quotient unstable strata}. 
			
			\bthm \label{thm can qnt unstab strata} \cite{Hoskins}
		 Let $G$ be a reductive group acting linearly on a projective scheme $X \subseteq \mathbb{P}^n$. Let $\|.\|$ be a Weyl-invariant inner product on $\Lie G$, and consider the associated HKKN stratification \[X=\bigsqcup\limits_{\b \in \mathscr{B}} S_\b.\]
			Suppose that for some $\b \in \mc{B} \setminus \{0\}$ such that we have condition (\ref{*}) as in Theorem \ref{Uhat thm}, i.e.\ that \[\stab_{U_\b}(x) = \{e\} \quad \forall x \in \Zb. \] Then the following holds:
			\bnu [label=\arabic*)] \item  The locus $\overbar{Y}_\b^{P_\beta-ss}$ as in Theorem \ref{Uhat thm} satisfies $\overbar{Y}_\b^{P_\beta-ss} \subseteq \Ybss$, and has a projective good quotient \[\overbar{Y}_\b^{P_\beta-ss} \rar \overbar{Y}_\b \git_{\mc{L}_{\beta, \eps}} \Pb \]
			\item We have, for $0<\eps \ll 1$, \[\overbar{Y}_\b \git_{\mc{L}_{\beta, \eps}} \Pb = \Proj \bigoplus_{i=0}^\infty H^0(\overbar{Y}_\beta,\mc{L}^{\otimes i}_{\beta,\eps}),\] where the quotient map is defined by the inclusion of the invariant subalgebra, which is finitely generated. The quotient is independent of $\eps$, taken to be strictly positive and sufficiently small.
		
			\item This quotient restricts to give a quasi-projective geometric quotient of the locus
			\[ \Yb^{ts} = \Yb^s \setminus U Z_{\beta}^s \rar (\Yb^s \setminus U Z_{\beta}^s)/P_\beta \subseteq \overbar{Y}_\b \git_{\mc{L}_{\beta, \eps}} \Pb  \]

			\item  Furthermore, suppose condition (\ref{ss=s R}) holds, i.e.\ that we have \[Z_\b^{\stab\b/\l_\b(\Gm)-s} = Z_\b^{\stab\b/\l_\b(\Gm)-ss}\]for the action of $\stab\b/\l_\beta(\Gm)$ with respect to the linearisation $\mc{L}_\b$. In this case, we have 
			$\overbar{Y}_\b^{P_\beta-ss}  =  Y^{ts}_\beta$ so that these geometric and projective quotients coincide. 
		
			\enu 
			\ethm
			
			\bpf
			
			We use Theorem \ref{Uhat thm} and the Hilbert-Mumford criterion as stated in Theorem \ref{thm HM crit}. We fix $\beta \in \mc{B} \setminus  \{0\}$, and aim to provide a quotient of $\Ybss$ by its action of $P_\beta$, which will be our $H$. Since we have a linearisation of the action of $G \acts X$, we have in particular an action of $P_\beta \leq G$ on the ambient projective space, and by the standard properties that an orbit closures is union of orbits and hence setwise invariant (see e.g.\cite{Brion} 1.11), we see that $\overbar{Y}_\beta$ is setwise invariant under this $P_\beta$ action. So $P_\beta$ acts on the projective variety $\overbar{Y}_\b$, which will play the role of $X$. 
			
		Up to conjugacy, there is a natural choice of distinguished cocharacter, namely an optimal cocharacter for $\Ybss$. To get a conjugacy class representative, we take $T$ a maximal torus of $G$, and intersect with $\Pb$ to get a maximal torus $T_\b$ of $\Pb$; that is, a maximal torus $T_\beta \leq \stab\beta$. Let $\l_\b :\Gm \rar P_\b$ be, amongst the $P_\b$-conjugacy class denoted by $\b$, the unique representative that lies in $T_\b$. This is the cocharacter we will use to grade $U_\b$. Indeed, since $P_\b = P(\l_\b)$ as defined in \S \ref{subsec HKKN}, all the weights of the adjoint action of $\l_\b$ on $\Lie P_\b$ must be non-negative. This establishes that the unipotent radical $U_\b \leq \Pb$ is graded by $\l_\b$.

			For the reductive part, observe that $\l_\b(\Gm)$ is by definition central in $\stab\beta$, so there is a character $\chi_\b$ of $\Pb$, given by the composition \[\chi_\b : \Pb \rar L_\b \rar \Gm.\] The first map in the composition is just the quotient map with kernel $U_\b$, and the second is obtained by extending to $\stab\b$ the character of the torus $T_\b$ dual to $\l_\b$ with respect to our chosen inner product; this is possible because $\l_\b$ is central in $L_\b$. We then define the linearisation $\mathscr{L}_{\b,\eps}$ to be $\mathscr{L}$ twisted by the character $(1+\eps)(-\chi_\b)$. By construction, for $0  <\eps \ll 1 $ this is adapted, since the minimal weight of $\l_\beta$ is $\beta$ and all other weights are strictly greater. In fact, by the arguments of \cite{Berczi2020} the choice of $\eps$ will not matter so long as it is taken sufficiently small and strictly positive. 
			
			With all of this in place, the theorem follows from the first two statements in Theorem \ref{thm HM crit}.

			\epf
			
			\brm
			Our construction of moduli spaces for unstable sheaves is, morally speaking, a matter of applying Theorem \ref{thm can qnt unstab strata}, but it not quite a straightforward application of this theorem, because we do not have condition (\ref{*}). However, for length two sheaves we will use Theorem \ref{Uhat for abelian U} instead, essentially proving a variant of Theorem \ref{thm can qnt unstab strata} that holds in this setting.
			\erm

			\subsection{Toy Example: Unordered points on $\pno$}
			\label{subsec pnts on pno}
			Now that the ideas have been introduced, the general situation of quotienting unstable strata is best clarified by working all of this out in detail for a simple example. Our description of the HKKN stratification here follows \cite{Kirwan1984}. Take $G \coloneqq \SL(2)$ and consider the natural action on \[X \coloneqq \text{Sym}^n(\pno) \cong \mathbb{P}^n\] induced by the action on $\pno$. We think of $X$ as the space of degree $n$ homogeneous polynomials $F(x,y)$, with the zeros of the polynomial giving the choice of points. This comes equipped with the action of $G$ induced by the dual action on $\pno$.  That is, $$(g \cdot F)(x,y) = F(g^{-1}\cdot(x,y)).$$ The linearisation is that induced from the standard $\mc{O}(1)$-linearisation on $\pno$: in representation-theoretic terms, this is $\sym^nV$, where $V$ is the two-dimensional irreducible representation. Our chosen maximal torus $T$ is identified with $\Gm$ via the isomorphism \[ T \rar \Gm\]
			$$\begin{pmatrix}
			t & 0 \\
			0 & t^{-1}
			\end{pmatrix} \longmapsto t.$$ This torus acts on $X$ diagonally with respect to the basis formed by monomials in $x$ and $y$, and hence the weights of $T$ are $$ \{ n-2i \in \ZZ  \mid i= 0,...n \} \subset \ZZ. $$
			
			We observe that if more than half of the points are at zero, there will be no negative weight, and if more that half are at infinity, there will be no positive weight. We therefore find, as is well known, that the first level of our stratification, the $G$-semistable locus, is: \[S_0 = X^{ss} = \{x \in X \mid \text{not more than half of the $n$ points coincide anywhere}  \}. \] We choose the positive Weyl chamber $\mathfrak{t}^+ \cong \RR_{\geq 0}$, and so the indexing set is:
			
			\[\mc{B} = \Big\{ 0 \Big\} \cup \Big\{ 2i-n \mid \frac{n}{2} < i \leq n \Big\}.  \] Immediately from the definition we compute that the parabolics and stabilisers are the same for every member of the indexing set, so that for all $i$ we have:
			
			\[P_{n-2i} = P \coloneqq \Big\{ \begin{pmatrix}
			a & b\\
			0 & a^{-1}
			\end{pmatrix} \mid a \in k^*, b \in k  \Big\} \] That is, the standard Borel subgroup in $\SL_2$. We also compute
			\[ \stab(n-2i) = T . \]  The $Z_{\beta}$'s each consist of a Stab$\beta$-invariant single point:
			$$Z_{n-2i} = \{ x \in X \mid \text{ Exactly }  i \text{ of the points are at $\infty$ and the rest are at $0$} \},$$ and we find: $$Y_{n-2i} = \{ x \in X \mid \text{ Exactly } i \text{ of the points are at $\infty$} \},$$ which we can check directly is $P$-invariant. 
			
			As a consequence, the reductive quotient of each stratum as described in \S\ref{hqotus} is just a point, which is clearly unsatisfactory. Let us see how using non-reductive GIT can improve on this.
			
			We consider the shifted linearisation $\mc{L}_\beta$ on a specific stratum, as used in \cite{Hoskins2011}. For a fixed $i$ we have the $T$-character $$\chi_{n-2i} : T \rar \Gm$$
			$$\begin{pmatrix}
			t & 0 \\
			0 & t^{-1}
			\end{pmatrix} \longmapsto t^{n-2i}$$
			and we define $\mc{L}_{n-2i}$ to be the linearisation of the $G$ action given by twisting $\mc{L}$ by $\chi_{n-2i}$. The condition (\ref{ss=s R}) holds trivially, since $\stab\b =\l_\b(\Gm)$. The (semi)stable locus of $Z_{\beta}$ with respect to $\mc{L}_{n-2i}$ is $$Z_{n-2i}^{s} = Z_{n-2i}^{ss}= Z_{n-2i} = \{\text{pt}\}.$$ Hence also, for each $\beta_i = (n-2i) \in \mc{B}$, we have $$Y_{n-2i}^{s} = Y_{n-2i}^{ss}=Y_{n-2i}.$$ Taking the projective closure yields: $$\overbar{Y}_{n-2i} = \{ x \in X \mid \text{At least } i \text{ of the points are at $\infty$} \} \cong \mathbb{P}^{n-i},$$ so, taking the $G$-sweep of the $Y_{n-2i}^{ss}$, the stratification is: $$S_{n-2i} = \{ x \in X \mid \text{Exactly } i \text{ of the points coincide}. \}$$
			
			Now we can apply the Theorem \ref{thm can qnt unstab strata} to find quotients, once we have verified the remaining remaining hypotheses. The subgroup $P$ fixes the point $\infty \in \pno$, so $P$ acts on $\overbar{Y}^{ss}_{n-2i}$ as stated above. Indeed, $\overbar{Y}^{ss}_{n-2i} \cong \mathbb{P}^{n-i}$ and the action is the natural one, and indeed it extends to an SL$(2)$ action on each $\mathbb{P}^{n-i}$. Now take $$ \begin{pmatrix}
			0 & \eta \\ 0 & 0
			\end{pmatrix} \in \text{Lie}(\Ga) \cong k$$ $$g_\lambda = \begin{pmatrix}
			\lambda & 0 \\ 0 & \lambda^{-1}
			\end{pmatrix}. $$ We see adjoint action of $\Gm$ on the Lie algebra of $\Ga$ is via $$g_\lambda \cdot X \cdot g_\lambda^{-1} = \begin{pmatrix}
			\lambda & 0 \\ 0 & \lambda^{-1}
			\end{pmatrix} \begin{pmatrix}
			0 & \eta  \\ 0 & 0
			\end{pmatrix}  \begin{pmatrix}
			\lambda & 0 \\ 0 & \lambda^{-1}
			\end{pmatrix}^{-1} = \begin{pmatrix}
			0 & \lambda^2 \eta  \\ 0 & 0
			\end{pmatrix}.$$ So the $\Gm$ in the Levi decomposition of $P$ acts on the Lie algebra Lie$(\Ga)$ of the $\Ga$ with sole weight $2 > 0$, in agreeement with the above. 
			
			Now we turn to weights of the $\Gm$ action from our linearisation. We began with the linearisation $\mc{L}$, with respect to which the weights on $\overbar{Y}_{n-2i}$ are  $ \{ (2i-n),... n \}$, and twist by the rational character associated to \[(1 +\eps)(n-2i),\] where $\eps >0$ is a small rational number. 
			
			It remain to verify the final hypothesis, to determine whether (\ref{*}) holds. For this it suffices to note that no point of $Z_{n-2i}$ is fixed under the translating action of the $\Ga$, unless $i=n$ and we are considering the case where all points coincide at $\infty$. This $Y_\beta$ corresponds to the $S_\beta$ stratum where all points coincide, and it is clear in this case that the quotient ought to be a point. At all other levels of the stratification, the hypotheses are verified and we can apply the Theorem \ref{thm can qnt unstab strata}.

			The theorem tells us that we have
			
			$$Y_{n-2i}^{ts}= (\overbar{Y}_{n-2i}) \setminus (\Ga \cdot  Z_{n-2i}).$$ Now we recognise the $\Ga$-sweep of $Z_{n-2i}$ as being the locus where $i$ points are at $\infty$ and the remaining points are affine and all coincide somewhere. So: $$Y_{n-2i}^{ts}=\overbar{Y}_{n-2i}^{P-ss}= \{ \text{Exactly $i$ points at $\infty$ and the remaining points don't all coincide} \}.$$ Since Theorem \ref{thm can qnt unstab strata} gives us a geometric quotient of this locus, we see that using non-reductive GIT gives a far better answer than the methods of \S \ref{hqotus}.

			\section{Moduli of Semistable Sheaves}\label{ss sh section}
			
			For the remainder of this article, we will be concerned with the moduli of coherent sheaves on a projective scheme. Our ultimate goal is to apply the results of \S \ref{sec moduli of unstable object} to quotient unstable strata in Simpson's reductive GIT construction of the moduli space of semistable sheaves, to obtain moduli spaces of unstable sheaves of fixed \HN type, in the case of \HN length 2. In this section we give a very brief review of semistable sheaves and Simpson's moduli construction \cite{Simpson1994}. The standard reference for this classical theory is \cite{Huybrechts2010}, and a very nice account of its historical development is given in \cite{Newstead1978}, \S5.1.
			
			Let $X$ be a projective scheme over $k$ with a fixed ample line bundle $\cO_X(1)$. For a coherent sheaf $F$ on $X$, we let $P(F,t) \in \QQ[t]$ denote the Hilbert polynomial of $F$ with respect to $\cO_X(1)$. As the Hilbert polynomial of such a sheaf is locally constant in flat families, from a moduli perspective it is natural to fix this discrete invariant. Hence, we fix $P \in \QQ[t] $, and from now on, unless specified otherwise, all coherent sheaves will have this Hilbert polynomial. 
			
			\subsection{Semistable sheaves}\label{moduli ss sheaves}  \label{sheaves lin} 
			We recall the definition of semistability for sheaves. We define a partial ordering $\preccurlyeq$ on $\QQ[t]$ by saying that $P \preccurlyeq Q$ whenever
			\[  \frac{P(n)}{P(m)} \leq \frac{Q(n)}{Q(m)} \] for $m \gg n \gg 0$. We define $\prec$ analogously by replacing $\leq$ with $<$. 
			
			\bdf \label{ss sheaf defn}
			
			We say $F$ is semistable  if for all non-zero proper subsheaves $F' \lneq F$, we have $P(F') \preccurlyeq P(F)$. If we have $P(F') \prec P(F)$ for all non-zero proper subsheaves, we say $F$ is stable.  
			
			\edf
			
			\brm The definition of semistability above is due to Rudakov \cite{Rudakov1997}, and is equivalent to the notion of Gieseker-semistability used by Simpson.  In particular, to be semistable it is necessary that a sheaf be pure; that is, the dimensions of the supports of all nonzero subsheaves must be equal. Hence, if $X$ is a smooth projective curve, the above is also equivalent to the notion of slope semistability for vector bundles on a curve. Our reason for using this formulation of semistability is that it allows us to assign Harder-Narasimhan filtrations to any coherent sheaf, not just the pure ones that would be allowed under Gieseker's definition.
			\erm
			
			Every semistable sheaf has a Jordan-H\"{o}lder filtration (\cite{Huybrechts2010} I\S1.5) that is, a filtration such that every successive subquotient is stable, and all subquotients have the same Hilbert polynomial as $F$. This filtration, unlike the \HN filtration we will consider later, is not canonical. However, its associated graded sheaf is determined up to isomorphism by $F$. 
			
			\bdf We say two semistable sheaves are \emph{S-equivalent} if the associated graded sheaves from their Jordan-H\"{o}lder filtrations are isomorphic.\edf
			
			\brm The moduli-theoretic motivation for this definition is as follows. If a semistable sheaf $F$ can be given as an extension \[0 \rar F_1 \rar F \rar F_2 \rar 0\] of two other sheaves, $F_1$, $F_2$ also semistable, then once can construct a family $\mathscr{F}$  of sheaves over $X \times \mathbb{A}^1$ - parameterised for example by a suitable line in $\Ext^1(F_2,F_1)$ - such that $\mathscr{F}_0 \cong F_1\oplus F_2$, and for $t \neq 0$ we have $\mathscr{F}_t \cong F$. This is an example of so-called jump phenomena, and shows us that any separated moduli space including all semistable sheaves must identify $F$ and $F_1\oplus F_2$. In particular, it shows that any such moduli space must identify all sheaves in the same S-equivalence class.
			\erm

			\subsection{The GIT problem} \label{subsec Simpson's GIT setup} By the Le Potier estimates (\cite{Simpson1994} Theorem 1.1), the family of semistable sheaves over $X$ with Hilbert polynomial $P$ is bounded. Hence, for $n \gg 0$, all such sheaves are $n$-regular in the sense of Castelnuovo-Mumford (\cite{Mumford1966} Lecture 14), i.e.\ we have $H^i(X,F(n-i)) = 0$ for all $i>0$. As observed by Mumford (\emph{loc. cit.}), $n$-regularity implies that
			
			\bit
			\item The sheaf $F(n)$ is generated by global sections, so that the evaluation map 
			\[ H^0(F(n)) \otimes \cO_X(-n) \ra F\]
			is surjective, by right exactness of tensor product;
			\item $H^i(F(n)) =0$ for $i>0$, so $\dim H^0(F(n)) = P(F,n)$.
			\eit

			This means that, if we fix a $k$-vector space $V$ of dimension $P(n)$, and are prepared to make a choice of isomorphism $H^0(F(n)) \simeq V$, we can instantiate any $n$-regular sheaf $F$ with Hilbert polynomial $P$ as the image of a quotient map $$q : V \otimes \cO_X(-n) \rar F.$$ Hence we can associate to any such sheaf a point in the Quot scheme 	\[ [q : V \otimes \cO_X(-n) \rar F] \in \quot_n \coloneqq \quot(V \ten \cO_X(-n),P)\] parameterising quotient sheaves of $V \ten \cO_X(-n)$ with Hilbert polynomial $P$. By a quotient sheaf we mean nothing more than an equivalence class of surjections $q: V \ten \cO_X(-n) \rar F$, where two such maps $q_i: V \ten \cO_X(-n) \rar F_i$, $i=1,2$, are regarded as equivalent, and hence correspond to the same point in the Quot scheme, if they have the same kernel; or, equivalently, if there exists an isomorphism of sheaves $\phi : F_1 \rar F_2$ making the below diagram commute:
			\begin{equation} 
			\begin{tikzcd} 
			V \ten \cO_X(-n) \arrow[d, "q_1"] \arrow[equal]{r} & V \ten \cO_X(-n) \arrow[d, "q_2"] \\ F_1 \arrow[r,"\phi"] & F_2.  
			\end{tikzcd} 
			\end{equation}
			It will be helpful to the reader to bear in mind that the Quot scheme is projective but not irreducible. The points associated to $n$-regular sheaves will all lie in the open subscheme $$Q = Q_n \subset \quot_n$$ defined to be that locus of quotient $q :V \ten \cO_X(-n) \rar F$ such that $H^0(q(n))$ is an isomorphism. 
			
			The point of $Q$ associated to a sheaf is not canonical: it depends on the choice of isomorphism $H^0(F(n)) \simeq V$ above. In other words, a point of the Quot scheme corresponds to a sheaf $F$, together with a choice of basis of $H^0(F(n))$. This gives rise to a natural action of $\GL(V)$ on $\Quot_n$, induced by the action on $V$, which preserves $Q$. For \[[q :V \ten \cO_X(-n) \rar F] \in \Quot_n\] and $g \in SL(V)$, we define $g \cdot [q]$ to be the quotient map given by the composition
			
			$$
			\begin{tikzcd} 
			V \ten \cO_X(-n) \arrow[r, "g^{-1} \ten \text{id}"] & V \ten \cO_X(-n) \arrow[r, "q"] & F.  
			\end{tikzcd}
			$$
			
			For later convenience, we note the following properties of this action on $Q$. 
			
			\blm \label{pp -  gl orbits = sheaves up to isom}
			(\cite{Huybrechts2010}, \S4.3)
			
			\begin{enumerate}[label=\arabic*)]
				\item Orbits of closed points of $Q$ under the $\GL(V)$ action correspond to isomorphism classes of sheaves. 
				\item For a closed point  $[q] \in Q$ corresponding to a globally generated sheaf $F$, there is an isomorphism $\stab_{\GL(V)}(q) \cong \Aut(F)$. \label{GL stab sheaves}
			\end{enumerate}
			\elm
			\bpf 
			The first statement is evident from the discussion above. For the second statement we note that any automorphism $\phi : F \rar F$ is determined by the map of vector spaces \[H^0(\phi(n)): H^0(F(n)) \rar H^0(F(n)).\] Conversely, the isomorphism $H^0(q(n)):V \rar H^0(F (n))$ coming from the point $q \in Q$ allows us to view an element of $g \in \GL(V)$ as an isomorphism $$H^0(q(n)) \circ g \circ H^0(q(n))^{-1} : H^0(F(n)) \rar H^0(F (n)).$$ By the definition of equivalence of points in the Quot scheme given by the diagram above, this isomorphism of twisted global sections induces an isomorphism of $F$ if and only if $g \in \stab(q)$. 
			\epf
			
			\brm \mbox{}
			\begin{enumerate}[label=\arabic*)]
				\item From the above proposition, it is clear that the centre $\mathcal{Z}(\GL(V))$ of $\GL(V)$ acts trivially, and so in our GIT set-up we will pass to considering an action of $SL(V)$. Up to replacing the linearisation with an appropriate tensor power, which we are certainly prepared to do, this is equivalent to considering the action of $\GL(V)/\mathcal{Z}(\GL(V)) = \PGL(V)$. See \cite{Huybrechts2010} 4.3 for more discussion of this.
				\item The above properties may be interpreted stack-theoretically, as saying that we have a presentation of $\cC oh_{X,P}$, the stack of coherent sheaves on $X$ with Hilbert Polynomial $P$, as an infinite increasing union of quotient stacks. This approach is taken, for example, in \cite{Hoskins2018}, \cite{Berczi2018}. For this it is necessary that we work with $\GL(V)$, since we must record all automorphisms of the sheaves, including the $\Gm$ common to all. 
				
			\end{enumerate} \erm

			The final ingredient necessary to turn this set-up into a reductive GIT problem is a linearisation. Take $m>n$ and write $H \coloneqq H^0(\mc{O}_X(m-n)),$ following \cite{Hoskins2012}. Now, we have a map $$ \quot_n \rar \text{Gr}(V \ten H,P(m))$$ $$[q : V \ten \mc{O}_X(-n) \rar F] \longmapsto [V\ten H \twoheadrightarrow H^0(F(m))].$$ to the Grassmannian of $P(m)$-dimensional vector space quotients $$V \ten H \rar W.$$  Furthermore, because Ker$(q)$ will be $m$-regular for $m$ large enough, if we take $m \gg n \gg 0$, this map is an embedding (\cite{Huybrechts2010}, 2.2.4). Following this up with the usual Pl\"{u}cker embedding of the Grassmannian $$ \text{Gr}(V \ten H,P(m)) \hrar \mathbb{P}\big( \bigwedge^{P(m)}(V \ten H)  ^* \big) \eqqcolon \PP^N$$ we obtain an equivariant embedding $i_{m,n}: \quot_n \hrar \PP^N$, which we may restrict to $Q$.  Denote by $\overbar{Q}$ the closure of $Q$ in this projective space. Since $\quot_n$ is projective ( see e.g.\ \cite{Fantechi2005}), we have $\overbar{Q} \subset \Quot_n$. Let us denote by $\mathscr{L}_{m,n}$ the very ample line bundle on $\Quot_n$ obtained from pulling back $\mc{O}_{\PP^N}(1)$ along this embedding.  \\~\\ We are now in the best setting for reductive GIT: we have a projective scheme $\overbar{Q}$ carrying an action of a reductive $G$, and the very ample line bundle corresponding to this embedding has a natural linearisation of the $G$ action. This is the set-up that Simpson uses to prove the following theorem.
			
			\bthm \cite{Simpson1994} \label{thm moduli of semistable sheaves} Let $X$ be a projective scheme, and fix a Hilbert Polynomial $P \in \QQ[t]$. Then
			
			\begin{enumerate}[label=\arabic*)]
				\item If $m \gg n \gg 0$, the semistable locus for the GIT problem above is   $$\overbar{Q}^{ss} = \{ [q : V \ten \mc{O}_X(-n) \rar F] \in Q \mid F \text{ is a semistable sheaf in the sense of \ref{ss sheaf defn}} \}.$$ Similarly, the GIT stable locus consists of those points of $Q$ corresponding to stable sheaves.
				\item The GIT quotient $M_{P,ss} \coloneqq \overbar{Q}\git_{\mathscr{L}_{m,n}} G$ is a coarse moduli space for $S$-equivalence classes of semistable coherent sheaves on $X$ with Hilbert polynomial $P$.
				
			\end{enumerate} 
			
			\ethm
			
			This moduli space, the natural successor to the moduli space of semistable vector bundles on a smooth curve, is by now a central object in algebraic geometry, and is very well studied from many different perspectives. From our point of view, however, it is just a starting point. Our eventual goal is to extend the classification of coherent sheaves on $X$ to include the unstable ones as well, so that Simpson's moduli space may be thought of as a special case: that of \HN length one.

			\section{Stratifications for Sheaves} \label{sheaves strats}

			Theorem \ref{thm moduli of semistable sheaves} gives a very satisfactory tool for understanding semistable sheaves and their moduli, but says nothing about sheaves that are \emph{unstable}, i.e.\ do not satisfy Definition  \ref{ss sheaf defn}. In understanding those, the following is the natural preliminary step.
			
			\subsection{Harder-Narasimhan Stratification}
			
			To any sheaf, one may associate a canonical filtration by subsheaves, first introduced by Harder and Narasimhan in \cite{Harder1975}. A proof of its existence and uniqueness in our setting may be found in \cite{Huybrechts2010}.
			
			\bdf \label{coprime HN defn} 
			Let $F$ be a sheaf on $X$. The \emph{Harder-Narasimhan filtration} (or \emph{HN filtration}) of $F$ is the unique filtration by subsheaves $$0 = F^0 \lneq F^1 \lneq \dots \lneq F^{\ell}  = F $$ such that the subquotients $F_i \coloneqq F^{i}/F^{i-1}$ are semistable with strictly decreasing Hilbert polynomials $P(F_1) \succ P(F_2) \succ \dots \succ P(F_{\ell})$. The \emph{associated graded} of $F$ is $\gr(F) \coloneqq \oplus_{i=1}^{\ell} F_i.$ We say that $F$ has \emph{\HN type} $\tau = (P(F_1),\dots,P(F_{\ell}))$, and \HN length $\ell$. \edf
			
			\bntt \emph{From now on, to prevent needless repetition, given a sheaf $F$ the terms in the \HN filtration of $F$ will always be denoted by $F^i$, and the subquotients will always be denoted by $F_i$.}
			\entt
			
			\brm By additivity of Hilbert Polynomials in exact sequences, if $F$ is of type $\tau = (P_1,\dots,P_\ell)$, then $\sum_i P_i = P$. \erm

			By results of Shatz and Nitsure (\cite{Shatz1977, Nitsure2011}), there is an upper semi-continuous \lq HN type function', which, given any family of sheaves $\mathscr{F} \rar S$, and a point $s \in S$, gives the \HN type of $\mathscr{F}_s$. Hence, defining \[\quot_n^\tau \coloneqq \{ [q: V \ten \cO_X(-n) \rar F] \in \quot_n \mid F\text{ has HN type $\tau$}\},\] one has a HN stratification\footnote{Sometimes also called the Shatz stratification.} into locally closed subschemes \[ \quot_n = \bigsqcup_{\tau} \quot^{\tau}_n.\] This in turn gives a stratification \[ Q = \bigsqcup_{\tau} Q_{\tau} \] defined in the same way. 
			
			\bdf We let $S_\tau$ denote the open subscheme of $Q_\tau$ consisting of those quotients $[q :V \ten \cO_X(-n) \rar F]$ such that $H^0(q(n))$ is an isomorphism.
			
			\edf 
			
			\brm
			All of the above can be translated into the language of stacks. It was shown in \cite{Hoskins2018} that there is a \HN type stratification on the stack, which may be thought of, in an appropriate sense, as a kind of colimit of the \HN stratifications of the Quot schemes, just as the stack itself is in some sense their colimit.
			\erm

			\subsection{Instability stratification for sheaves} \label{instab strat sheaves}
			In addition to the stratification by \HN type one also has, according to the general principles described in \S\ref{subsec HKKN}, an HKKN instability stratification of the Quot scheme, \[ Q = \bigsqcup_{\beta \in \cB(n,m)} S_\beta \] coming from Simpson's GIT set-up.\footnote{This stratification depends in a non-trivial way upon the parameters $m$ and $n$, so strictly speaking the first question one should ask is whether there is any meaningful sense in which these stratifications can be said to stabilise. Fortunately, as we shall see, there is.} The natural thing to do, then, is to compare the two. Indeed, since they both measure the failure of a sheaf to be semistable, one might hope that they would coincide. The relationship between these stratifications was studied by Hoskins and Kirwan. First, in \cite{Hoskins2012}, it was shown that for a pure HN type\footnote{That is, a HN type of a pure sheaf.} $\tau$  and for $m\gg n\gg 0$ depending on $\tau$, there is a conjugacy class $\beta_{n,m}(\tau) $ of rational cocharacters of $G$ such that the HN stratum $Q^\tau$ is contained (set theoretically) in the instability stratum $S_{\beta_{n,m}}(\tau)$. In \cite{Hoskins2018}, this result was extended to a scheme theoretic statement for all HN types of coherent sheaves. Theorem \ref{thm comp strat} summarises these results. First, give the definition of the HKKN index $\beta_{n,m}(\tau)$ corresponding to a HN type $\tau$.
			
			\begin{defn}\label{def hess index from HN type} \cite{Hoskins2018} 
				For a \HN type $\tau = (P_1, \dots , P_l)$, we define a rational $\SL_{P(n)}$-cocharacter $\l_{\beta(n,m,\tau)}: \Gm \SL_{P(n)}$, such that $\l_{\beta(n,m,\tau)}(t)$ is diagonal with entries \[(t^{\b_1},\dots t^{\b_1},t^{\b_2},\dots t^{\b_2},\dots, t^{\b_\ell},\dots, t^{\b_\ell})\] where each \[\beta_i:= \frac{P(m)}{P(n)} - \frac{P_i(m)}{P_i(n)}\] is repeated $P_i(n)$ times. We write $\beta(n,m,\tau)$ for the conjugacy class of this cocharacter under conjugation by its associated parabolic $P_{\l_{\beta(n,m,\tau)}}$, as defined in \S\ref{subsec HKKN}.

			\end{defn}
			
			\brm
			Instead of working with rational cocharacters -that is, cocharacters with rational weights- we could instead define $\l_{\beta(n,m,\tau)}(t)$ to be the unique indivisible cocharacter obtained as a positive rational multiple of the cocharacter defined above. \erm 
			
			\brm
			Because two values will not always uniquely specify a polynomial, even for $m \gg n \gg 0$ the map $\tau \mapsto \beta(n,m,\tau)$ will not be injective if $\dim X > 1$. However, in \cite{Hoskins2018} Hoskins and Kirwan have shown that there is an asymptotic HKKN stratification on the stack of coherent sheaves, and this coincides with the stratification by HN type.
			\erm 
			
			A cocharacter $\lambda:\Gm \rar \SL(V)$ is the same thing as the data of its weight filtration \[\{0\} \subseteq \dots \subseteq V^{i} \subseteq V^{i-1}\subseteq \dots \subseteq V\] where $V^{i}$ denotes the span of those vectors on which $\lambda$ acts with weight at least $i$, so that only finitely many of the $ V_i :=\faktor{V^i}{V^{i-1}}$ are non-zero. Thus, given a point \[[q :V \ten \cO_X(-n) \rar F] \in \Quot_n,\] such a weight filtration induces a filtration of $F$, by defining \[F^i \coloneqq q(V^{i}\otimes\mc{O}_X(-n)).\] Limits under cocharacters are easy to compute.
			
			\blm \label{rmk limits in the quot scheme} (\cite{Huybrechts2010} Lemma 4.4.3) We have \[\lim_{t \rar 0} \l(t)[q] = [\oplus_i q_i :V \ten \cO_X(-n) \rar \oplus_i F_i]\in\overbar{Q}\] where for each $i$, \[q_i: V_i \ten \cO_X(-n) \rar F_i\] is the surjection induced from $q$. In words: $F$ is mapped to the associated graded of the filtration determined by $\l$. 
			\elm 
			
			The origin of the cocharacter $\l_{\b(n,m,\tau)}$, then, is not so mysterious. In looking for optimal cocharacters for $F$, a sensible initial guess is that this filtration should be $F$'s \HN filtration.  As in \cite{Hoskins2011}, one then obtains $\l_{\b(n,m,\tau)}$, by maximising (with Lagrange multipliers) the normalised Hilbert-Mumford function over all weight filtrations that induce the \HN filtration of $F$.

			The fixed point locus $\overbar{Q}^{\l_\b}$ thus consists of all quotients \[[\oplus_i q_i :V \ten \cO_X(-n) \rar \oplus_i F_i]\in\overbar{Q}\] that split with respect to the weight filtration of $\l_\b$. Because the Hilbert polynomials of the summand image sheaves must be locally constant, there is a component $F$ of $\overbar{Q}^{\l_\b}$ where for each $i$, the sheaf $F_i$ has Hilbert polynomial $P_i$. And so, writing, \[ \Quot_{n,i} = \Quot(V_i\otimes \cO_X(-n),P_i) \] one has:
			\[F \cong  \Quot_{n,1} \times \dots \times  \Quot_{n,\ell}.\]
			
			As above, we let $Q_i \subset  \Quot_{n,i}$ be the open locus of pure dimension $d$ sheaves such that $H^0(q_i(n))$ is an isomorphism. This allows one to define:
			\[Z_\tau = \{\oplus_i q_i \in F \mid q_i \in Q_i \}.\]
			Thus $Z_\tau$ is open in $F$, but its closure is a proper subscheme. Now let $Q^{ss} \subset Q_{i}$ be the open locus consisting of those sheaves which are in addition semistable, and define
			\[Z^{ss}_\tau = \{\oplus_i q_i \in F \mid q_i \in Q^{ss}_i \}.\]
			
			Note that all sheaves in $Z^{ss}_\tau$ are of \HN type $\tau$, and isomorphic to their \HN graded. We now have the notation we need to state the results of \cite{Hoskins2012} \cite{Hoskins2018}.

			\begin{thm}[\cite{Hoskins2012,Hoskins2018}]\label{thm comp strat}
				Let $\tau=(P_1, \dots , P_l)$ be a HN type for sheaves on $X$ with Hilbert polynomial $P$; then, for $m\gg > n\gg 0$, the following statements hold for $\beta:=\beta(n,m,\tau)$.
				\begin{enumerate}[label=\roman*)]
					\item Isomorphism classes of sheaves of \HN type $\tau$ are in bijection with $G$-orbits of points in $S_\tau$.
					\item The conjugacy class $\beta \in \cB$ is an instability index for the $G$ action on $Q$ with respect to $\mathscr{L}_{n,m}$ and $S_\tau$ is a closed subscheme of the corresponding instability stratum $S_{\beta}$.
					\item We have $Z_\tau^{ss} \subset Z_\beta^{ss}$ as a $\stab\b$-invariant closed subscheme.
					\item If $Y_\tau^{ss} := p_\beta^{-1}(Z_\tau^{ss})$, then $S_\tau = G \cdot Y_\tau^{ss} \cong G \times^{P_\b} Y_\tau^{ss}.$
					\item The subcheme $Y_\tau^{ss}$ is precisely the locus of points  $[q:V \ten \cO_X(-n) \rar F] \in \Quot$ such that $H^0(q(n))$ is an isomorphism, $F$ has \HN type $\tau$, and the weight filtration of $\l_\b$ induces the \HN filtration of $F$. The map $p_\b : Y_\tau^{ss} \rar Z_\tau^{ss}$ takes these sheaves to their associated graded.
				\end{enumerate}
			\end{thm}

			We can give explicit GIT characterisations of the subschemes $Y_\tau^{ss}$ and $Z_\tau^{ss}$, in the manner usual for instability stratifications. A detailed account is in \cite{Hoskins2011}; we give an overview here to fix notation for later.

			\brm \label{rmk on structure of Pbeta} We recall from Remark \ref{rmk filtration parabolics} that $P_\beta = P_{\l_{\beta}}$ is the subgroup of $\SL(V)$ consisting of elements that respect the weight filtration of $\l_{\beta(n,m,\tau)}$, i.e.\ block upper triangular matrices. Thus $P_{\l_{\beta}}$ splits as a semidirect product of a reductive part, $\stab \beta$, which consists of block diagonal matrices, and its unipotent radical $U_\beta$, which consists of block upper triangular matrices with identity matrices on the diagonal blocks. Thus one has 
			
			\[\stab\b \cong  \left( \prod_i \GL(V_i) \right) \cap \SL(V), \]
			
			and the action on $Z_\tau$ is via each $\GL(V_i)$ acting on $Q_i$ in the way described at \S\ref{subsec Simpson's GIT setup}. \erm 
			To analyse the GIT semistability with respect to $\stab\b$, \cite{Hoskins2012} considers the subgroup \[\prod_{i=1}^{\ell} \SL_{P_i} \leq \stab\b.\] This group acts on \[Z_\tau \subset Q_1\times \dots \times Q_\ell\] diagonally via the natural action of $\SL_{P_i}$ on $Q_i$, discussed in \S\ref{subsec Simpson's GIT setup}. Furthermore, the linearisation induced for this action by the restriction of the linearisation $\mc{L} = \mc{L}_{m,n}$, is nothing but \[ \mc{L}_\pi :=
			\pi_1^*\mc{L}_1 \otimes \dots \otimes \pi_{\ell}^* \mc{L}_\ell \] where $\mc{L}_i$ is the linearisation given at \S\ref{subsec Simpson's GIT setup}, and $\pi_i: \prod_j Q_j \rar Q_i$ is the projection. 
			
			If $Z(\stab\b) \cong \Gm^{\ell-1}$ denotes the centre of $\stab\b$, then one has a surjection with finite kernel \[\prod_{i=1}^{\ell} \SL_{P_i}\times Z(\stab \b) \rar \stab \b. \] Twisting the linearisation by the character $\b$ ensures that the centre $Z(\stab \b)$ acts trivially on the canonical linearisation $\mc{L}_\b$ of Definition \ref{canonical lin defn}. Hence, since GIT (semi)stability is unaffected by finite groups, one has:
			
			\[Z_\tau^{\stab \beta-ss}(\mc{L}_\b) = Z_\tau^{\prod \SL(V_i)-ss}(\mc{L}_\pi),\]
			i.e.\ the semistable loci for $\stab \beta$ with respect to $\mc{L}_\b$, and $\prod_i \SL_{P_i}$ with respect to $\mc{L}_\pi$, coincide. From the form of the linearisation we have, for any cocharacter \[\l = (\l_1,\dots \l_\ell): \Gm \rar \prod_i \SL_{P_i}\] and any point $q = (q_1,\dots, q_{\ell}) \in Z_\tau$, that \[\mu_{\mc{L}_{\pi}}(q,\l) = \sum_{i=1}^{\ell} \mu_{\mc{L}_i}(q_i,\l_i).\] From this one readily deduces: 		     
			\bpp \cite{Hoskins2012, Hoskins2018} \label{prop Ztauss is semistable locus for Stab beta}	The subscheme $Z_\tau^{ss} \subset Z_\tau$ is exactly the semistable locus for the $\stab\b$ action on $Z_\tau$ with respect to the linearisation $\mc{L}_\beta$. \epp	
			
			Thus, following \S \ref{hqotus} it is possible to obtain categorical quotients of $Y_\tau^{ss}$ using only reductive GIT, by taking the $P_\beta$-quotient of $Y_\tau^{ss}$ using the canonical linearisation $\mc{L}_\b$. 
			However, \S \ref{hqotus}  tells us that if we do this we will learn nothing new: we identify each sheaf with its accociated \HN graded, and hence the resultant \lq moduli spaces of unstable sheaves' are isomorphic to products of moduli spaces of semistable sheaves. We now turn to non-reductive GIT and the methods of \S \ref{subsec nr qnts of unstable strata}, to see if this situation can be improved.

			\section{Moduli of $\tau$-stable Sheaves of \HN length 2}\label{main section}
			
			Having set up all of the machinery for constructing moduli spaces of unstable objects via non-reductive GIT, we now turn to our main application: sheaves of Harder-Narasimhan length two. Throughout we fix a \HN type $\tau = (P_1,P_2)$, of length $\ell =2$. Fix $m \gg n \gg 0$ and let $\beta = \beta(n,m,\tau)$ be the associated HKKN index.	
			
					\subsection{On $\tau$-stable sheaves} \label{subsec on tau-stable sheaves}
			\bdf \label{defn coprime} We say that a \HN type $\tau = (P_1, \dots , P_\ell)$ is coprime\footnote{The reason for this choice of term is that coprimality generalises the situation of semistability and stability coinciding for vector bundles on a curve, which happens when the rank and degree, two positive integers, are coprime in the usual sense.} if, for all $i = 1, \dots \ell$, all semistable sheaves on $X$ of Hilbert Polynomial $P_i$ are stable. \edf
			
			We now define the sheaves for which we will construct moduli spaces. 
			 	\bdf \label{defn tau-stable sheaf} We say that a sheaf $F$ of HN length 2 is $\tau$-stable if it is of \HN type $\tau$, we have $F 
			 \ncong \gr(F)$, and the HN-subquotients $F_i$ in its associated graded $\gr(F) \cong \oplus_{i=1}^{\ell} F_i$, are stable. \edf

			Later on we will give a GIT interpretation of such sheaves in the HKKN picture for Simpson's construction. This will allows us to prove our main result, Theorem \ref{intro main thm}, which we repeat here for the reader's convenience.

			\bthm \label{main thm} Let $(X,\mc{O}_X(1))$ be a projective scheme with an ample line bundle, and let $\tau =(P_1,P_2)$ be a Harder-Narasimhan type of length 2. Then 
			\bnu 
			\item There is a quasi-projective moduli space $M_{\tau,d}$ for those $\tau$-stable (Definition \ref{defn tau-stable sheaf}) coherent sheaves such that $\dim\End(E) = d$. 
			\item This moduli space has a canonical projective completion $M_{\tau,d} \subset \overbar{M}_{\tau,d}$, constructed via the nonreductive partial desingularisation procedure of \S \ref{blow up discussion}, such that the map \[M_{\tau,d} \rar M^{ss}_{P_1}\times M^{ss}_{P_2},\] to the product of two moduli spaces of semistable sheaves, which takes a sheaf to the associated graded of its HN filtration, extends to $\overbar{M}_{\tau,d}$.
			\item Furthermore, if $\tau$ is coprime in the sense of Definition \ref{defn coprime}, then $M_{\tau,d}$ is the moduli space of indecomposable sheaves of type $\tau$ such that $\dim\End(E) = d$. 
			\enu 
			\ethm

			We pause to make the following elementary observation.
			\blm  \label{lem tau stables are indec}
			A sheaf of \HN length 2 is indecomposable if and only if $\gr F \ncong F$. In particular, $\tau$-stable sheaves are indecomposable. If $\tau$ is coprime, the converse holds as well.
			\elm 
			\bpf
			Of the first statement, one implication is clear. For the other, recall that every stable sheaf $E$ is simple, (see Corollary 1.2.8 of \cite{Huybrechts2010}) i.e.\ it satisfies $\End(E) \cong k$, so in particular is indecomposable. Hence if there were a splitting $F \cong F^\prime \oplus F^{\prime \prime}$, then $F^1$ would be entirely contained within one factor, say $F^\prime$. But then $ F_2 \cong F^{\prime}/F^1\oplus F^{\prime}$, contradicting the simplicity, and hence the stability, of $F_2$. 
			
			For the last statement, any indecomposable sheaf of type $\tau$ must have non-split \HN filtration, and the condition that the subquotients $F_i$ are stable holds trivially because $\tau$ is coprime.
			\epf 
			
			As might be expected, there is a GIT characterisation of $\tau$-stability, analogously to \S\ref{instab strat sheaves}.
			
			\bdf  Write $Z_\tau^{s}$ for the subscheme of $Z_\tau^{ss}$ consisting of HN-graded sheaves whose HN-summands are all stable. \edf 
			
			The following is proved in the same way as Proposition \ref{prop Ztauss is semistable locus for Stab beta}.
			\bpp 
			The locus $Z_\tau^{s}$ is precisely the stable locus for the action of $\overbar{R}_\b =  \faktor{\stab \b}{ \l_\b(\Gm)}$ on $Z_\tau$, with respect to the linearisation $\mc{L}_\beta$. \label{prop ztau stable locus}  
			\epp 	
			\bpf
			Recall the discussion preceding Proposition \ref{prop Ztauss is semistable locus for Stab beta}, and in particular the surjection \[\prod_{i=1}^{\ell} \SL_{P_i}\times Z(\stab \b) \rar \stab \b. \] Since we are in \HN length 2, we have $\l_\b(\Gm) \cong Z(\stab \b)$, so we get a surjection \[\prod_i \SL_{P_i} \rar \faktor{\stab \b}{\l_\b(\Gm)} = \overbar{R}_\b \] which also has finite kernel. Since GIT stability is also unaffected by finite group actions, we have an equality of stable loci, \[Z_\tau^{\overbar{R}_\beta-s}(\mc{L}_\b) = Z_\tau^{\prod_i\SL_{P_i}-s}(\mc{L}_\b).\] Hence recalling the Hilbert-Mumford criterion for stability from \S\ref{recap of reductive GIT}, it suffices to show that for all nontrivial \[\l= (\l_1,\dots \l_\ell): \Gm \rar \prod_i \SL_{P_i}\] and any point $q = (q_1,\dots, q_{\ell}) \in Z_\tau$, that $\mu_{\mc{L}_{\b}}(q,\l) >0.$ By the decomposition \[\mu_{\mc{L}_{\pi}}(q,\l)=\sum_{i=1}^{\ell} \mu_{\mc{L}_i}(q_i,\l_i),\] this will hold if and only if $\mu_{\mc{L}_i}(q_i,\l_i)>0$ for each $i$. And this is the case if and only if each $q_i$ is a stable sheaf- that is, if and only if $q\in Z_\tau^s$.

			\epf
			This allows us to deduce a GIT characterisation of the sheaves for which we will construct a moduli space.	
			
			\bcr \label{cor GIT characterisation of tau-stable}
			If we define $Y^{s}_\tau := p_\beta^{-1}(Z_\tau^{s})$, then $P_\beta$-orbits in $Y^s_\tau \setminus U_\b Z_\tau^s$ are in one-to-one correspondence with isomorphism classes of $\tau$-stable sheaves.
			\ecr 			
			
			\brm A word of warning is in order here. It might seem that the above proposition is almost a tautology, and should therefore be true in for any HN length. However, this is not the case, and the reason is somewhat subtle. Consider a point $q \in Z_\tau^s$, where now $\tau$ has length $\ell$. The corresponding HN-graded sheaf $F \cong \gr(F) = \oplus_{i=1}^{\ell} F_i$ of \HN length $\ell$ has a copy of $\Gm^{\ell}$ in its automorphism group (and hence, by Lemma \ref{pp -  gl orbits = sheaves up to isom} its $\GL(V)$-stabiliser in the Quot scheme) which comes from scaling each factor. 
			
			Passing to $\SL(V)$ removes the diagonal $\Gm$ that is common to all sheaves, and so we are left with a copy of $\Gm^{\ell-1}$ inside the $\stab \b$-stabiliser of such a sheaf. Upon quotienting by $\lambda_{\beta}(\Gm)$, this reduces to a $\Gm^{\ell-2}$-dimensional torus inside the $\overbar{R}_\beta$-stabiliser of $q \in Z_\tau^{s}$, meaning that such points cannot be $\overbar{R}_\beta$-stable if $\ell >2$. Put another way, in the language of the convex hull formulation of the Hilbert-Mumford criterion, the origin is contained in the interior of the weight polytope when the interior is taken with respect induced topology on the polytope itself, but not when the interior is taken with respect to the topology on the codimension-one subspace of $\mathfrak{t}^*$ corresponding to the Lie algebra of a maximal torus of $\overbar{R}_\beta$.  This is one reason for the alternative construction pursued in \cite{Hoskins}, which constructs moduli spaces for certain sheaves of \HN length $\ell >2$.

			\erm

			\subsection{Automorphism groups of sheaves of \HN length 2}
			\label{subsec automorphism gps of unstable sheaves}
			We are now almost ready to prove Theorem \ref{main thm}. 
			
			Following \S \ref{sec Uhat}, we  will need to understand the unipotent radicals of stabilisers of points in the Quot scheme corresponding to sheaves \HN length 2. By Lemma \ref{pp -  gl orbits = sheaves up to isom}, the $\GL$-stabiliser of such a point is isomorphic to the automorphism group of the sheaf, and so we collect some observations about automorphism groups of unstable sheaves. We begin with the following, doubtless well known, which shows that there is a natural parabolic structure on these automorphism groups.
			
			\blm \label{lem sheaf auts preserve hn filt} Let $F$ be a coherent sheaf on $X$. Then all automorphisms of $F$ preserve the \HN filtration of $F$.
			\elm
			\bpf 
			This follows directly from the uniqueness of the \HN filtration. 
			\epf

			\brm 
			 The above lemma may be seen as an instance of a more general phenomenon occurring in the setting of moduli spaces of unstable objects. In view of the GIT construction in \S \ref{ss sh section} and \S\ref{sheaves strats}, the statement for automorphisms follows from \cite{Kempf1978}, so long as we are careful to use the action of $\GL(V)$ rather than $\SL(V)$ for the GIT setup, i.e.\ to have the right stack-theoretic parameter space.
			\erm

			Next we need some notation for versions in $\GL(V)$ of the groups we have encountered.
			
			\bntt
			We write $\widetilde{P}_\b \cong P_\b \times Z(\GL(V))$ for the subgroup of $\GL(V)$ generated by $P_\b$ and the diagonal $\Gm$. Thus, $\widetilde{P}_\b \cong \U_\b \rtimes \widetilde{\stab}\b$, where $\widetilde{\stab}\b$ is the group generated by $\stab\b$ and the central $\Gm$ in $\GL(V)$. In other words, we let $\widetilde{P}_\b$ be the parabolic group defined by the weight filtration of $\l_\b(\Gm)$, taken in $\GL(V)$.
			\entt 
			
			\brm 
			Note that, since the centre of $\GL(V)$ acts trivially on the Quot scheme, the orbits for these groups are the same as for their analogues in $\SL(V)$. The only difference is the action of $Z(\Gm)$ on the linearisation, which is by the determinantal character of $\GL(V)$. Thus, we could just as well set up the GIT problem using these groups, so long as we twist the linearisation we are using by that character. Indeed, as we have seen, from a stack-theoretic perspective this contruction is more natural, since it gives the correct automorphism groups.
			\erm
			
			Now fix a sheaf $F$ of type $\tau = (P_1,P_2)$, represented by a choice of point in the Quot scheme $[q: V\otimes \mc{O}_X(-n) \rar F] \in Y_\tau^{ss}$. This choice allows us to make free use of the identification \[H^0(q(n)): V  \rar H^0(F(n)),\] and gives us a faithful representation $\Aut(F) \hrar \GL(V)$, restricting to an isomorphism \[\stab_{\GL(V)}(q) \cong \Aut(F)\] as in Lemma \ref{pp -  gl orbits = sheaves up to isom}.

				We write \[\Aut(F) \cong \Aut_U(F) \rtimes \Aut_R(F),\] where $\Aut_U(F) \trianglelefteq \Aut(F)$ is the unipotent radical, and $\Aut_R(F)$ is a chosen Levi factor. By the above Lemma \ref{lem sheaf auts preserve hn filt}, in fact $\Aut(F) \leq \widetilde{P}_\b$, from which it follows that $\Aut_U(F) = \Aut(F) \cap U_\b$, and we can choose $\Aut_R(F) =  \Aut(F) \cap \widetilde{\stab} \beta$.

			\bpp (\emph{cf. \cite{BrambilaPaz2013}})
			Let $F$ be a sheaf of \HN length $\ell=2$.
			Then we have \[\Aut_U(F) \cong \id_F + \Hom(F_2,F_1),\] as groups, with composition being the operation on the right.  
			Furthermore, if $F$ is $\tau$-stable then we have $\End(F) \cong k \oplus \Hom(F_2,F_1).$

			\epp 
			
			\bpf

			Using the representation of $\Aut(F)$ above, from $\Aut_U(F) =  \Aut(F) \cap U_\b$ we see \[\Aut_U(F) \subseteq \{\id_F + \phi \in \End(F) \mid \phi \in \Hom(F,F_1)\}.\] In fact, since such an element must have $F_1$ in its kernel by its block form, we see that $\phi \in \Hom(F_2,F_1)$, where we use the natural injection $\Hom(F_2,F_1) \hrar \Hom(F,F_1)$ which comes from applying $\Hom(-,F_1)$ to the short exact sequence $0 \rar F_1 \rar F \rar F_2\rar 0.$ Conversely any $\phi \in \Hom(F_2,F_1)$ yields an element $\id_F + \phi \in \Aut_U(F)$ in this manner.
			
		Now suppose that $F$ is $\tau$-stable By Lemma \ref{lem tau stables are indec} it is then indecomposable, and by \cite{Atiyah1956} this means that we have \[\End(F) \cong k\oplus \Nil(E)\] where $\Nil(E)$ is the nilpotent endomorphisms of $E$. There is an injective map \[\Hom(E_2,E_1) \rar \Nil(E)\] sending $\phi \in \Hom(E_2,E_1)$ to the composition
			\[\begin{tikzcd}
			E \arrow[twoheadrightarrow]{r}  & E_2 \arrow{r}{\phi}  & E_1 \arrow[hookrightarrow]{r} &E.
			\end{tikzcd}\]
			Now for surjectivity. Take $\alpha \in \Nil(E)$.  We observe that since all automorphisms of $F$ preserve the HN filtration, by considering $\id_F +\alpha$, so must $\alpha$. Thus we get $\alpha|_{E^1}: E^1 \rar E^1$, which must be zero since $E^1$ is simple and $\alpha$ is nilpotent. By the same reasoning, the induced map $\overbar{\alpha}: E_2 \rar E_2$ must be zero, which shows that $\alpha \in \Hom(E_2,E_1)$.  
				\epf
		
Putting the pieces together, we find the following important corollary,  which greatly simplifies the blow-up process, and gives a sheaf-theoretic intrinsic description of the resultant refined \HN type.	
			\bcr \label{cor pbeta preserves ustab dim}
			\bnu \item  For any $q \in Y_\tau^{ss}$ corresponding to a sheaf $F$ of type $\tau$, we have \[\dim \stab_{U_\b}(q) = \dim \Hom(F_2,F_1)\] and hence \[\dim \stab_{U_\b}(q) = \dim \stab_{U_\b}(p_\b(q)),\] so condition \ref{$p$ preserves U stab dim} holds.
			
			\item If in addition $F$ is $\tau$-stable, we have \[\dim \End(F) = \dim \stab_{U_\b}(q) +1.\]  
			
			\enu 
			 \ecr  
			\brm \label{rmk higher length stab interp} The failure of this Corollary for higher length \HN types, and its consequences for the blow-up process described at \S \ref{blow up discussion}, is an important reason why in this article we restrict our attention to length two. The analogous statement in higher length $\ell$ is that there is a filtration 
						\[0 = \End^\ell(F) \subsetneq \dots \subsetneq \End^1(F) \subsetneq \End(F)\] where \[\End^k(F) := \{ \phi \in \End(F) \mid \phi(F^{i}) \subseteq F^{i-k} \quad \forall i=1,\dots \ell\},\] 		
			and we have $\Aut_U(F) \cong \id_F + \End^{1}(F)$. For the graded sheaves it follows, using the standard property (\cite{Huybrechts2010} Prop 1.2.7) that there are no homomorphisms \lq forwards' in the HN filtration, that we have \[\Aut(\gr F) \cong \left( \prod_{i=1}^{\ell} \Aut(F_i) \right) \ltimes \bigoplus_{i>j}\Hom(E_i,E_j),\] and hence that \[\Aut_U(\gr F) \cong \id_{\gr F} + \bigoplus_{i>j}\Hom(E_i,E_j).\] 
			Similarly to the above, we can get an injective map \[\End^{1}(F) \hrar \bigoplus_{i>j}\Hom(E_i,E_j).\] However, there is no reason to suppose that this map will be surjective, meaning that the analogue of Condition \ref{$p$ preserves U stab dim} need not hold.

\erm

			\subsection{Proof of the main theorem}

			Now that all of the necessary ideas are in place, we are ready to prove the main result. 
						The basic idea of the proof is this: beginning with a suitable subset \[Y^{\prime} \subset \overbar{Y_{\tau}}\] of the projective closure of $Y_{\tau}$, one performs the blow-up process described above in \S \ref{blow up discussion}, as in Theorem \ref{Uhat for abelian U}. This corresponds to refining the instability stratification, and hence to refining the \HN stratification of the Quot scheme. This results in a space $\pi: \widehat{Y}^{\prime} \rar Y^{\prime}$ satisfying the condition \ref{*_d}. An open subset of $\widehat{Y}^{\prime}$ therefore has a projective categorical quotient by the $P_\b$ action, and these quotients will be the projective completions mentioned in Theorem \ref{main thm}. To prove the theorem, we must show that certain loci in this projective completion can be identified with the right moduli spaces of sheaves. This will be done via the isomorphism over the complement of the exceptional divisor given by restricting the blow-down map, as in Corollary \ref{cr uh4}. The crucial point is that, as a result of the first part of Corollary \ref{cor pbeta preserves ustab dim}, we are in the good situation discussed in \S \ref{blow up discussion}, and so the locus we blow up will not at any stage contain the minimal weight space. This allows us to explicitly determine a locus that remains stable and untouched throughout the blow-up process, whose quotient is the quasi-projective moduli space of the theorem. Before moving on to the proof of the theorem, let us make a few remarks.
						
						\brm \mbox{}
						
						\begin{enumerate}[label=\arabic*)]
							\item It may happen for a given value of $\tau$ and $d$ that $M_{\tau, d}$ is empty: for example, if $X$ is a smooth projective curve and we choose $\tau$ to be such that $\Ext^1(F_2,F_1)$ is forced to vanish for degree reasons. 
							\item As observed above, the situation is particularly good if the \HN type in question is coprime. Similarly to the reductive-semistable case, this is a consequence of the fact that, for these HN types, semistability and stability (in the sense of GIT) coincide for the action of $\overbar{L_\beta}$. And, again as in the reductive case \cite{Kirwan1985} if this does not happen one can perform a sequence of blow-ups according to the dimensions of reductive stabilisers, to obtain a space on which this condition is satisfied, yielding partial desingularisations of the projective completions in the theorem above, and hence an alternative compactification of $M_{\tau,d}$.

						\end{enumerate}

						\erm

						\begin{proof}[Proof of Theorem \ref{main thm}]
							
First, we note that the unipotent radical is abelian, by its description at Remark \ref{rmk on structure of Pbeta}.  Thus, we blow up according to the dimensions of stabilisers in $U$. Now, let $d_{\min}$ and $d_{\max}$ be respectively the minimal and maximal values of $\dim \stab_U(y)$ across $ y \in Y_\tau$.

By Corollary \ref{cor pbeta preserves ustab dim}, since we consider only $q \in Y_\tau^{ss}$, corresponding to a point in the Quot scheme $[q : V \ten \mc{O}(-n) \rar \cF]$, this is equivalent to blowing up according to the dimension of the group $\Hom(\cF_2,\cF_1)$. 
							\\~\\
\noindent \underline{\emph{Case: $d = d_{\min}$}} \mbox{} \\~\\ Consider the locus $$(Y^{s}_\tau)^{d_{\min}} \coloneqq \{ y \in Y^{s}_\tau \mid \dim \stab_U(y) = d_{\min}  \}.$$ Then by Proposition \label{prop ztau stable locus} and Corollary \ref{cor pbeta preserves ustab dim}, orbits under the $P_\b$ action on this space correspond to isomorphism classes of sheaves $\cF$ of \HN type $\tau$, with $\dim \Hom(\cF_2,\cF_1) = d_{\min}$, and $F_1$, $F_2$ stable sheaves. Furthermore, the sheaves in $$P_\beta (Z_\tau^s)^{d_{\min}} = U_\b(Z^{s}_\tau)^{d_{\min}} = U_\b Z^{s}_\tau \cap (Y^{s}_\tau)^{d_{\min}} $$ are precisely those sheaves which are in addition \HN split. Hence, by Corollary \ref{cor GIT characterisation of tau-stable}, to construct the moduli space $M_{\tau,d}$ it suffices to exhibit a geometric quotient of the locus $(Y^{s}_\tau)^{d_{\min}} \setminus U_\beta(Z^{ss}_\tau)^{d_{\min}}$ by its $P_\b$ action.  \\~

Taking the projective closure $\overbar{Y}_\tau$ in the ambient projective space, we first blow up $\overbar{Y_\tau}$ along the closure of the subset $$C_{d_{\max}}(Y_\tau) = \{ y \in Y_\tau \mid \dim \stab_{\Uh}(y) = d_{max} \} = U\{y \in Z_{\min} \mid \dim\stab_U(y) = d_{\max}-1\} $$ This results in a space $\overbar{Y_\tau}^{(1)}$ with an action of $P_\b$ lifting that on $\overbar{Y_\tau}$. Since the centre of the blow-up did not contain $Z_\tau$, the minimal weight space for $\lambda_\beta$ acting on $\overbar{Y_\tau}^{(1)}$, which we denote by $Z_{\tau}^{(1)}$, is the proper transform of the minimal weight space $Z_\tau$ of $\overbar{Y_\tau}$, as in \S\ref{blow up discussion} case $(1)$. Let us simplify notation by denoting by $Y_\tau^{(1)}$ the open subset $$(\overbar{Y_{\tau}}^{(1)})^0_{\min} \subset \overbar{Y_{\tau}}^{(1)}$$ consisting of the points flowing down to $Z_{\tau}^{(1)}$ under $\lambda_\tau$ as $t\rar 0$. Finally, we linearise the $P_\b$ action on $\overbar{Y_\tau}^{(1)}$ by replacing the linearisation $\mathscr{L}_{\tau}^{per}$ with its perturbation by a small multiple of the exceptional divisor of the blow-up. Because all of these blow-ups may be performed in the ambient space, we may assume that at each stage of the blow-up process the exceptional divisor is just the projectivised normal bundle to the centre of the blow-up, and the corresponding weights are just certain small perturbations of the weights for $Y$, in the way described in Remark \ref{blow up wts Gm}.
							\\~\\
							This process strictly reduces the dimensions of the stabilisers in $U_\b$ over the centre of the blow-up (see \cite{Berczi2020} and also \cite{Hoskins} \S 2.3.7), thus guaranteeing that the dimensions $$\dim \stab_{U_\b}(y),\text{ for }y  \in Z_\tau^{(1)}$$ are all strictly less than $d_{\max}$. We now blow up $\overbar{Y_{\tau}}^{(1)}$ along the closure of the locus of points where the $\Uh_\beta$-stabiliser is maximal, and we repeat this process until it terminates with a space $ \pi : \widehat{Y}_\tau \rar \overbar{Y_\tau}$. This space carries a $P_{\tau}$ action lifting that on $\overbar{Y_\tau}$, and this action has a well-adapted linearisation $\mathscr{\widehat{L}}^{per}_{\tau}$ obtained as described, by perturbing at each stage of the blow-up process with a small multiple of the exceptional divisor. Let us denote the  minimal $\lambda_\tau$ weight space in $\widehat{Y}_\tau$ by $\widehat{Z}_\tau$. This has basin of attraction $(\widehat{Y}_{\tau})^0_{\min}$. Then we have arrived at a situation where $$\dim \stab_{U_\b}(y) = d_{\min}\text{ } \forall\text{ }y\in (\widehat{Y}_{\tau})^0_{\min}. $$ In other words, condition \ref{*_d} is satisfied, so by Theorem \ref{Uhat for abelian U}  there is a projective categorical quotient of the locus $\widehat{Y}_\tau^{P_\b-ss}$, given by composing a geometric quotient by $\Uh_\beta$ with a categorical quotient by the residual action of the reductive group $(\stab\beta)/\lambda_{\beta}(\Gm)$. We define $$ \overbar{M}_{\tau,d_{\min}} \coloneqq \widehat{Y}_\tau \git_{\mathscr{\widehat{L}}^{per}_{\tau}} P_\b .$$ In addition, the complement $U_\b\widehat{Z}_{\min}$ also has a categorical quotient by the $P_\b$ action. Indeed, the geometric quotient of $U_\b\widehat{Z}_{\min}$ by $\Uh$ is of course $\widehat{Z}_{\min}$, and the open subset $\widehat{Z}_{\min}^{\stab \beta-ss}$ has a categorical quotient given by the usual GIT quotient with respect to the canonical linearisation, as described in \S\ref{hqotus}. Hence we observe that $\widehat{Z}_{\min} \git_{\mathscr{L}_\beta} \stab\beta$ is a categorical quotient of  $U\widehat{Z}_{\min}$. Moreover, the blow-down map $ \pi : \widehat{Y}_\tau \rar \overbar{Y_\tau}$ restricts to give a blow down \[\widehat{Z}_{\min} \rar Z_\tau, \] again because the new minimal weight space is the proper transform of the old one. Moreover, under this map, the image of the semistable locus for $\stab\beta \acts \widehat{Z}_{\min}$ with respect to $\widehat{\mathscr{L}}_\beta$ is contained in $Z_\tau^{ss}$ - that is, the semistable locus for $\stab\beta \acts Z_\tau$ with respect to $\mathscr{L}_\beta$. This follows by the same argument as in 
							\cite{Kirwan1985} because the linearisation $\widehat{\mathscr{L}}_\beta$ may be chosen to be an arbitrarily small peturbation of the pullback $\pi^*\mathscr{L}_\beta$ by the exceptional divisor.

							 Since the blow-up is equivariant, we thus get a map of their quotients \[ \widehat{Z}_{\min} \git_{\widehat{\mathscr{L}}_\beta} \stab\beta \rar Z_\tau\git_{\mathscr{L}_{\beta}}\stab\b \cong M_{P_1}^{ss} \times M_{P_2}^{ss},\] where the latter is a product of Simpson moduli spaces of semistable sheaves. It remains to identify a suitable quasi-projective subscheme of $\overbar{M}_{\tau,d_{\min}}$ with the desired moduli space of sheaves.

							Let $E$ be the exceptional locus of the map $\pi : \widehat{Y}_\tau \rar \overbar{Y_\tau}$; then $\pi$ gives an identification $$\widehat{Y}_\tau \setminus E \cong \overbar{Y_\tau} \setminus \pi(E).$$ The locus $(Y^{ss}_\tau)^{d_{\min}}$ is untouched by the blow-up, and hence this isomorphism restricts to an identification of $(Y^{ss}_\tau)^{d_{\min}}$ with a subset of $\widehat{Y}_\tau \setminus E$.  And, since at every stage $\pi_{i} : \overbar{Y_{\tau}}^{(i+1)} \rar \overbar{Y_{\tau}}^{(i)}$ of the blow-up process the centre of the blow-up did not contain the minimal $\lambda_\tau$ weight space of $\overbar{Y_{\tau}}^{(i)}$, the blow down map $\pi : \widehat{Y}_\tau \rar \overbar{Y_\tau}$ also restricts to an isomorphism $\widehat{Z}_\tau \setminus E \cong (Z_\tau)^{d_{\min}}$. Futhermore, all of these identifications are $P_\b$-equivariant, since $\pi$ is $P_\b$ equivariant.
							
							~\\ Proceeding analogously to \S \ref{subsec on tau-stable sheaves},  let $\widehat{Z}^{s}_\tau \subset \widehat{Z}_\tau$ denote the stable locus for the action of $\stab\beta$ on $\widehat{Z}_\tau$ with respect to the canonical linearisation $\widehat{\mathscr{L}}_\beta$ for $\widehat{Y}_\tau$. Then let $\widehat{Y}^{s}_\tau$ be the inverse image $\widehat{p}_\b^{-1}(\widehat{Z}^{ss}_\tau)$ under the flow $$\widehat{p}_{\b} :  (\widehat{Y}_\tau)^0_{\min} \rar \widehat{Z}_\tau$$ given by taking $t \rar 0$ in the $\lambda_\tau$ action. By the Hilbert-Mumford criterion given in Theorem \ref{thm HM crit} applied to the blow-up space $\widehat{Y}_\tau$, one has a quasi-projective geometric $P_\b$-quotient \[\widehat{Y}_\tau^{ts} \rar \widehat{Y}_\tau^{ts}/P_\b \subseteq \overbar{M}_{\tau,d_{\min}}. \]

							  Then, by the classical Hilbert-Mumford criterion for reductive groups, we may identify $\widehat{Z}^{s}_\tau \setminus E$ with $(Z^{s}_{\tau})^{d_{\min}}$, via $\pi$. It follows that we may identify $(Y^{s}_{\tau})^{d_{\min}}$ with $\widehat{Y}^{s}_{\tau} \setminus E$, and hence we obtain via the blow-down map a $P_\b$-equivariant isomorphism \[\widehat{Y}^{ts}_\tau \setminus E \cong (Y_{\tau}^{ts})^{d_{\min}}. \] Taken together, this tells us we can restrict the $P_\b$-quotient map $$\widehat{Y}_\tau^{P_\b -ss} \rar \overbar{M}_{\tau,d_{\min}} $$ to get a quasi-projective geometric $P_\b$-quotient:  $$M_{\tau,d_{\min}} \coloneqq  (Y^{ts}_\tau)^{d_{\min}}/P_\b.$$ This proves the first part of the theorem in the case $d=d_{\min}$. 
							  
							  For the second part, observe that that the map \[\pi \circ \widehat{p}_\tau: (\widehat{Y)}^0_{\min} \rar \widehat{Z}_\tau \rar Z_\tau\] is $P_\b$-equivariant, and that we have $\widehat{Y}_\tau^{ts} \subset (\widehat{Y}_\tau)^0_{\min}$. Hence the above map, may be restricted to  $\widehat{Y}_\tau^{ts}$, and then induces a map of quotients by the universal property of categorical quotients:  $\overbar{M}_{\tau,d_{\min}} \rar M_{P_1}^{ss} \times M_{P_2}^{ss},$ which extends the map taking a sheaf to its associated graded.

							  For the final part of the theorem, we use the last part of Lemma \ref{lem tau stables are indec}.
							  \\

							\noindent \underline{\emph{Case: $d > d_{min}$}} \mbox{}
							\\~\\
							Supposing inductively that the theorem holds for some $d-1$, we return to our original $\overbar{Y}_{\tau}$. By upper semi-continuity of stabiliser dimension, the subset $$\overbar{Y^{d}_{\tau}} \coloneqq \overbar{Y}_{\tau} \setminus \{ y \in \overbar{Y_\tau} \mid \dim \stab_U(y) < d  \}$$ is closed in $\overbar{Y}_{\tau}$. Furthermore, condition \ref{$p$ preserves U stab dim} ensures that, whenever the locus \[(Y_\tau^s)^{d}) =  Y_\tau^d \cap Y_\tau^{s}\] is non-empty, the locus $(Z_\tau^s)^{d}= Z_\tau^d \cap Z_\tau^{s}$ will also be non-empty. We then simply repeat the process described above, replacing $\overbar{Y}_{\tau}$ by $\overbar{Y^{d}_{\tau}}$ and $d_{\min}$  by $d$.

							\epf

							In particular, as a special case we recover a result obtained in \cite{BrambilaPaz2013} by different means. It is interesting to note that the same discrete invariant appears there as when one approaches the problem using non-reductive GIT.
							
							\bcr \cite{BrambilaPaz2013}
							Let $X$ be a nonsingular projective curve, and $\tau = (P_1,P_2)$ a pure HN type corresponding to coherent sheaves of rank 2. Then the space $M_{\tau,d}$ above is a quasi-projective moduli space of indecomposable vector bundles $F$ of \HN type $\tau$ with $\dim\End(\cF) = d+1$. 
							\ecr

							\bibliography{ell=2}
							
							\bibliographystyle{alpha}
							
						\end{document}